\documentclass[aap,MSNbibl,seceqn,citesort,dvips]{arximspdf}
\usepackage{dcolumn}

\doi{10.1214/10-AAP758}
\volume{21}
\issue{6}
\pubyear{2011}
\firstpage{2315}
\lastpage{2342}

\makeatletter

\newtheorem{theorem}{Theorem}
\newtheorem{lem}{Lemma}
\newproclaim{exa}{Example}

\newcolumntype{d}[1]{D{.}{.}{#1}}

\newcommand{\var}{\operatorname{Var}}
\newcommand{\cF}{\mathcal{F}}
\newcommand{\cX}{\mathcal{X}}
\newcommand{\bX}{\mathbf{X}}
\newcommand{\bY}{\mathbf{Y}}
\newcommand{\wht}{\widehat}
\newcommand{\wtd}{\widetilde}
\newcommand{\wQ}{\widetilde Q}
\newcommand{\wX}{\widetilde X}
\newcommand{\wY}{\widetilde Y}
\newcommand{\wbY}{\widetilde{\mathbf Y}}
\newcommand{\rD}{\mathrm{D}}
\newcommand{\rB}{\mathrm{B}}
\newcommand{\rR}{\mathrm{R}}

\makeatother

\begin{document}
\begin{frontmatter}

\title{A sequential Monte Carlo approach to
computing tail probabilities in stochastic
models}
\runtitle{Sequential Monte Carlo for tail probabilities}

\begin{aug}
\author[A]{\fnms{Hock Peng} \snm{Chan}\corref{}\thanksref{t1}\ead[label=e1]{stachp@nus.edu.sg}} and
\author[B]{\fnms{Tze Leung} \snm{Lai}\thanksref{t2}\ead[label=e2]{lait@stat.stanford.edu}}
\runauthor{H. P. Chan and T. L. Lai}
\affiliation{National University of Singapore and Stanford University}
\address[A]{Department of Statistics\\
\quad and Applied Probability\\
National University of Singapore\\
6 Science Drive 2 \\
Singapore 117546\\
\printead{e1}}
\address[B]{Department of Statistics\\
390 Serra Mall\\
Stanford, California 94305\\
USA\\
\printead{e2}}
\end{aug}

\thankstext{t1}{Suported by the National University of Singapore Grant
R-155-000-090-112.}
\thankstext{t2}{Supported by the NSF Grant DMS-08-05879.}

\received{\smonth{2} \syear{2009}}
\revised{\smonth{12} \syear{2010}}

%
\begin{abstract}
Sequential Monte Carlo methods which involve sequential importance
sampling and resampling are shown to provide a versatile approach to
computing probabilities of rare events. By making use of martingale
representations of the sequential Monte Carlo estimators, we show how
resampling weights can be chosen to yield logarithmically efficient
Monte Carlo estimates of large deviation probabilities for
multidimensional Markov random walks.
\end{abstract}

%

\end{frontmatter}

\section{Introduction}

In complex stochastic models, it
is often difficult to evaluate probabilities of events of interest
analytically and Monte Carlo methods provide a practical
alternative. When an event $A$ occurs with a~small probability (e.g.,
$10^{-4}$), generating 100 events would require a very large number
of events (e.g., 1 million) for direct Monte Carlo computation of
$P(A)$. To circumvent this difficulty one can use importance
sampling instead of direct Monte Carlo changing the measure $P$ to
$Q$ under which $A$ is no longer a rare event and evaluating $P(A) =
E_Q(L {\mathbf1}_A)$ by $m^{-1} \sum_{i=1}^m L_i {\mathbf1}_{A_i}$, where
$(L_1,{\mathbf1}_{A_1}),\ldots,(L_m,{\mathbf1}_{A_m})$ are $m$ independent
samples drawn from the distribution $Q$, with $L_i$ being a
realization of the likelihood ratio statistic $L:=dP/dQ$, which is
the importance weight. While large deviations theory has provided
important clues for the choice of $Q$ for Monte Carlo evaluation of
exceedance probabilities, it has also been demonstrated that
importance sampling measures that are consistent with large
deviations can perform much worse than direct Monte Carlo (see
Glasserman and Wang \cite{GW97}). Chan and Lai \cite{CL07} have recently
resolved this problem by showing that certain mixtures of
exponentially twisted measures are asymptotically optimal for
importance sampling. For complex stochastic models, however, there
are implementation difficulties in using these asymptotically
optimal importance sampling measures. Herein we introduce a \textit
{sequential importance sampling and resampling} (SISR) procedure to
attain a weaker form of asymptotic optimality, namely, logarithmic
efficiency; the definitions of \textit{asymptotic optimality} and
\textit{logarithmic efficiency} are given in Section~\ref{sec3}.

Instead of applying directly the asymptotically optimal importance
sampling measure $Q$ that is difficult to sample from, SISR
generates $m$ sequential samples from a more tractable importance
sampling measure $\wQ$ and resamples at every stage $t$ the $m$
sequential sample paths, yielding a~modified sample path after
resampling. The objective is to approximate the target measure $Q$
by the weighted empirical measure defined by the resampling weights.
Details are given in Section \ref{sec2} for general resampling weights (not
necessarily those
associated with the asymptotically optimal resampling measure).
Section \ref{sec4} illustrates the SISR method for Monte Carlo
computation of
exceedance probabilities in a variety of applications which include
boundary crossing probabilities of generalized likelihood ratio
statistics and tail probabilities of Markov random walks. These
applications demonstrate the versatility of the SISR method and the
relative ease of its implementation.

Our SISR procedure to compute probabilities of rare events is
closely related to (a) the \textit{interacting particle systems} (IPS)
approach introduced by Del Moral and Garnier \cite{DG05} to compute tail
probabilities of the form $P \{ V(X_t)\,{\geq}\,a \}$ for a possibly
nonhomogeneous Markov chain $\{ X_t \}$ and (b)~the \textit{dynamic
importance sampling} method introduced by Dupuis and Wang \cite
{DW05,DW07} to
compute $P \{ S_n/n \in A \}$, where $S_n=\sum_{t=1}^n g(X_t)$ and
$\{ X_n \}$ is a uniformly recurrent Markov chain with stationary
distribution $\pi$ such that $\int g(x) \,d \pi(x) \notin A$. Both
IPS and dynamic importance sampling generate the $X_i$ sequentially.
Dynamic importance sampling uses an adaptive change of measures based
on the simulated
paths up to each time $t \leq n$. A recent method closely related to dynamic
importance sampling is sequential state-dependent change of measures
introduced by Blanchet and Glynn \cite{BG08} for Monte Carlo evaluation
of tail probabilities of the maximum of heavy-tailed random walks.
The IPS approach uses ``mutation'' to sample~$\wX^{(i)}_{t+1}$
(conditional on the
$X_1^{(i)},\ldots,X_t^{(i)}$ already generated) from the original
measure $P$ and then uses ``selection'' to draw $m$ i.i.d.\vspace*{1pt} samples from
$\{
(X_1^{(i)},\ldots,X_t^{(i)}, \wX_{t+1}^{(i)})\dvtx1 \leq i \leq m \}$
according to a Boltzmann--Gibbs particle measure.
The theory of IPS in \cite{DG05} focuses on tail probabilities of
$V(X_t)$ for fixed
$t$ as described in Section \ref{sec2} rather than large deviation probabilities
of $g(S_n/n)$ for
large $n$ as considered in Section \ref{sec3}.
Our SISR procedure is motivated by rare events of the general form $\{
\bX_n \in
\Gamma\}$ that involves the entire sample path $\bX_n = (X_1,\ldots,
X_n)$ and includes $\{ V(X_n) \geq a \}$ and $\{ S_n/n \in A \}$
considered by Del Moral and Garnier, Dupuis and Wang as special cases.
The sequential importance sampling component of SISR
uses an easily implementable approximation $\wQ$ of $Q$; in many
cases it simply uses $\wQ=P$. Thus, it is quite different from
dynamic importance sampling even though both yield logarithmically
efficient Monte Carlo
estimates of $P \{ S_n/n \in A \}$.

\section{Sequential importance sampling
and resampling (SISR) and martingale representations}\label{sec2}

The events in this section are
assumed to belong to the $\sigma$-field generated by $n$ random
variables $Y_1,\ldots,Y_n$ on a probability space $(\Omega,\mathcal
{F},P)$. Let $\bY_t=(Y_1,\ldots,Y_t)$ for $1 \leq t \leq n$. For
direct Monte Carlo computation of $\alpha:=P \{ \bY_n \in\Gamma
\}$, i.i.d. random vectors $\bY_n^{(1)},\ldots,\bY_n^{(m)}$ are
generated from $P$ and $\alpha$ is estimated by
%
%
\begin{equation} \label{2.1}
\wht\alpha_{\mathrm{D}} = m^{-1} \sum_{i=1}^m {\mathbf1}_{\{ \bY_n^{(i)}
\in\Gamma\}}.
\end{equation}
The estimate $\wht\alpha_{\rD}$ is unbiased and its variance is
$\alpha(1-\alpha)/m$ which can be consistently estimated by
%
%
\begin{equation} \label{2.2}
\wht\sigma_{\rD}^2 := \wht\alpha_{\rD} (1-\wht\alpha_{\rD})/m.
\end{equation}

In most stochastic models of practical interest, the $Y_t$ are
either independent or are specified by the conditional densities
$p_t(\cdot|\bY_{t-1})$ of $Y_t$ given~$\bY_{t-1}$, with respect to
some measure $\nu$. Direct Monte Carlo computation of $P \{ \bY_n
\in\Gamma\}$, therefore, involves $Y_1^{(i)},\ldots, Y_n^{(i)}$ that are
generated sequentially from~the\-se conditional
densities for $1\leq i\leq m$. In contrast, SISR first generates~%
$m$~inde\-pendent random variables $\wtd Y_t^{(1)}, \ldots,
\wtd Y_t^{(m)}$ at stage $t$, with~$\wtd Y_t^{(i)}$ having density
function $\widetilde q_t(\cdot|\bY_{t-1}^{(i)})$ to form $\wtd
\bY_t^{(i)}=(\bY_{t-1}^{(i)}, \wY_t^{(i)})$ and then uses resampling
weights of the form $w_t(\wtd\bY_t^{(i)})/ \sum_{j=1}^m
w_t(\wbY_t^{(j)})$ to draw $m$ independent sample paths
$\bY_t^{(j)}$, $1 \leq j \leq m$, from $\{ \wtd\bY_t^{(i)}, 1 \leq
i \leq m \}$. Here $\widetilde q_t$ are conditional density functions with
respect to $\nu$ such that $\widetilde q_t > 0$ whenever $p_t > 0$; one
particular choice is $\widetilde q_t=p_t$. In Section \ref{sec3}, we
show how the
weights~$w_t$ can be chosen to obtain logarithmically efficient SISR
estimates of rare event probabilities.

The preceding SISR procedure uses \textit{bootstrap resampling} that
chooses~i.i.d. sample paths from a weighted empirical measure of \mbox{$\{
\wtd\bY_t^{(i)}, 1 \leq i \leq m \}$}. It is, therefore, similar to
the selection step of the IPS approach that chooses i.i.d.
``path-particles'' from some weighted empirical particle measure (see~%
\cite{DG05}). The Monte Carlo estimate of $\alpha$
using SISR with bootstrap resampling is
%
%
\begin{equation} \label{2.3}
\wht\alpha_{\rB} = m^{-1} \sum_{i=1}^m L\bigl(\wbY_n^{(i)}\bigr)
h_{n-1}\bigl(\bY_{n-1}^{(i)}\bigr) {\mathbf1}_{\{ \tilde\bY_n^{(i)} \in\Gamma
\}},
\end{equation}
where $h_0 \equiv1$ and
%
%
\begin{eqnarray} \label{2.4a}
L({\mathbf y}_n)&=&\prod_{t=1}^n \frac{p_t(y_t|{\mathbf y}_{t-1})}
{\widetilde q_t(y_t|{\mathbf y}_{t-1})},\qquad h_k({\mathbf y}_k)= \prod_{t=1}^k
\frac{\bar
w_t}{w_t({\mathbf y}_t)},\nonumber\\[-6pt]\\[-6pt]
\bar w_t &=& \frac{1}{m} \sum_{i=1}^m
w_t\bigl(\wbY_t^{(i)}\bigr).\nonumber
\end{eqnarray}

Chan and Lai \cite{CL08} have recently developed a general theory of
sequential Monte Carlo filters in hidden Markov models by using a
representation similar to the right-hand side of (\ref{2.3})
for these filters. The method of their analysis can be applied
to analyze $m(\wht\alpha_{\mathrm B}-\alpha)$, decomposing it into
a sum of $(2n-1)m$ terms so that the summands form a martingale difference
sequence. Let $E^*$ denote expectation under the probability
measure $\wtd Q$ from which the $\wbY^{(i)}_t$ and $\bY_t^{(i)}$ are
drawn and define for $1 \leq t < n$,
%
%
\begin{equation} \label{2.5}
f_t({\mathbf y}_t)=E^*\bigl[L(\bY_n) {\mathbf1}_{\{ \bY_n \in\Gamma\}}|
\bY_t={\mathbf y}_t\bigr]=L({\mathbf y}_t) P(\bY_n \in\Gamma|\bY_t =
{\mathbf y}_t),
\end{equation}
setting $f_0 \equiv\alpha$ and $f_n(\wbY_n)=L(\wbY_n) {\mathbf1}_{\{
\tilde\bY_n \in\Gamma\}}$.
An important\vspace*{1pt} ingredient in the analysis is the ``ancestral origin''
$a_t^{(i)}$ of $\bY_t^{(i)}$. Specifically, recall that the ``first
generation'' of the $m$ particles consists of
$\wY_1^{(1)},\ldots,\wY_1^{(m)}$ (before resampling) and set
$a_t^{(i)}=j$ if the first component of $\bY_t^{(i)}$ is
$\wY_1^{(j)}$. Let $\#_k^{(i)}$ denote the number of copies of
$\wbY_k^{(i)}$ generated from $\{ \wbY_k^{(1)},\ldots,\wbY_k^{(m)}
\}$ to form the $m$ particles in the $k$th generation and let
$w_k^{(i)}=w_k(\wbY_k^{(i)})/\sum_{j=1}^m w_k(\wbY_k^{(j)})$. Then
it follows from (\ref{2.4a}) and simple algebra that for $1 \leq i
\leq m$,
\begin{eqnarray*} 
mw_t^{(i)} & = & h_{t-1}\bigl(\bY_{t-1}^{(i)}\bigr)/ h_t\bigl(\wbY_t^{(i)}\bigr), \\[3pt]
\sum_{i\dvtx a_t^{(i)}=j} f_t\bigl(\bY_t^{(i)}\bigr) h_t\bigl(\bY_t^{(i)}\bigr) & = &
\sum_{i\dvtx
a_{t-1}^{(i)}=j} \#_t^{(i)} f_t\bigl(\wbY_t^{(i)}\bigr)
h_t\bigl(\wbY_t^{(i)}\bigr),\\[-15pt]
\end{eqnarray*}
\begin{eqnarray*} 
&& \sum_{t=1}^n \sum_{i\dvtx a_{t-1}^{(i)}=j}
\bigl[f_t\bigl(\wbY_t^{(i)}\bigr)- f_{t-1}\bigl(\bY_{t-1}^{(i)}\bigr)\bigr]
h_{t-1}\bigl(\bY_{t-1}^{(i)}\bigr) \\[3pt]
&&\quad{} + \sum_{t=2}^n
\sum_{i\dvtx a_{t-2}^{(i)}=j} \bigl(\#_{t-1}^{(i)}-mw_{t-1}^{(i)}\bigr) f_{t-1}
\bigl(\wbY_{t-1}^{(i)}\bigr) h_{t-1}\bigl(\wbY_{t-1}^{(i)}\bigr) \\[3pt]
&&\qquad = \sum_{i\dvtx
a_{n-1}^{(i)}=j} f_n\bigl(\wbY_n^{(i)}\bigr) h_{n-1}\bigl(\bY_{n-1}^{(i)}\bigr) -\alpha,
\end{eqnarray*}
recalling that $f_0 \equiv\alpha$, $h_0 \equiv
1$ and defining $a_0^{(i)}=i$. Let
%
%
\begin{eqnarray} \label{2.10}
\varepsilon_{2t-1}^{(j)} &=&
\sum_{i\dvtx a_{t-1}^{(i)}=j} \bigl[f_t\bigl(\wbY_t^{(i)}\bigr)-
f_{t-1}\bigl(\bY_{t-1}^{(i)}\bigr)\bigr]h_{t-1}\bigl(\bY_{t-1}^{(i)}\bigr)
\quad\mbox{for } 1
\leq t \leq n, \nonumber\hspace*{-35pt}\\[-8pt]\\[-8pt]
\varepsilon_{2t}^{(j)} &=&
\sum_{i\dvtx a_{t-1}^{(i)}=j} \bigl(\#_t^{(i)}-mw_t^{(i)}\bigr) \bigl[
f_t\bigl(\wbY_t^{(i)}\bigr)h_t(\wbY_t^{(i)}) -\alpha\bigr] \quad\mbox{for } 1 \leq t
\leq n-1. \nonumber\hspace*{-35pt}
\end{eqnarray}
Then for each fixed $j$, $\{ \varepsilon_t^{(j)}, 1 \leq t
\leq2n-1 \}$ is a martingale difference sequence with respect to
the filtration $\{ \cF_t, 1 \leq t \leq2n-1 \}$ defined below and
%
%
\begin{equation} \label{2.9}
m(\wht\alpha_{\rB} - \alpha) = \sum_{j=1}^m \bigl(\varepsilon_1^{(j)} +
\cdots+ \varepsilon_{2n-1}^{(j)}\bigr).
\end{equation}
The martingale representation (\ref{2.9}) that involves the ancestral origins
of the genealogical particles is useful for
estimating the standard error of $\wht\alpha_{\mathrm B}$,
as shown by Chan and Lai \cite{CL08} who have also introduced the
$\sigma$-fields
%
%
\begin{eqnarray} \label{2.11}
\cF_{2t-1} & = &
\sigma\bigl( \bigl\{ \wtd Y_1^{(i)}\dvtx1 \leq i \leq m \bigr\} \nonumber\\
&&\hphantom{\sigma\bigl(}
{}\cup\bigl\{ \bigl(\bY_s^{(i)},
\wbY_{s+1}^{(i)},a_s^{(i)}\bigr)\dvtx1 \leq s < t, 1
\leq i \leq m \bigr\} \bigr), \\
\cF_{2t} & = & \sigma\bigl( \cF_{2t-1} \cup\bigl\{ \bigl(\bY_t^{(i)},a_t^{(i)}\bigr)\dvtx
1 \leq i \leq m \bigr\}\bigr) \nonumber
\end{eqnarray}
with respect to which (\ref{2.10}) forms a martingale difference sequence.

Since\vspace*{1pt} $f_n(\wbY_n^{(i)})=L(\wbY_n^{(i)}) {\mathbf1}_{\{ \tilde
\bY_n^{(i)} \in\Gamma
\}}$ and $\sum_{i=1}^m (\#_t^{(i)}-mw_t^{(i)})=0$ for $1 \leq t \leq n-1$,
summing (\ref{2.10}) over $t$ and $j$ yields (\ref{2.9}). Without
tracing their ancestral origins, we can also use the
successive generations of the $m$ particles to form martingale
differences directly.
Specifically, in analogy with (\ref{2.10}), define for
$i=1,\ldots,m$,
%
%
\begin{eqnarray} \label{list}
Z_{2t-1}^{(i)} &=& \bigl[f_t\bigl(\wbY_t^{(i)}\bigr)-f_{t-1}\bigl(\bY_{t-1}^{(i)}\bigr)\bigr]
h_{t-1}\bigl(\bY_{t-1}^{(i)}\bigr) \quad\mbox{for } 1 \leq t \leq n,
\nonumber\hspace*{-35pt}\\[-8pt]\\[-8pt]
Z_{2t}^{(i)} &=& f_t\bigl(\bY_t^{(i)}\bigr) h_t\bigl(\bY_t^{(i)}\bigr)-\sum_{j=1}^m w_t^{(j)}
f_t\bigl(\wbY_t^{(j)}\bigr) h_t\bigl(\wbY_t^{(j)}\bigr) \quad\mbox{for } 1 \leq t
\leq n-1. \nonumber\hspace*{-35pt}
\end{eqnarray}
As noted by Chan and Lai \cite{CL08}, $\{ (Z_t^{(1)}, \ldots, Z_t^{(m)}),
1 \leq t \leq2n-1 \}$ is a martingale difference sequence with
respect to the filtration $\{ \cF_t, 1 \leq t \leq2n-1 \}$ and $Z_t^{(1)},
\ldots, Z_t^{(m)}$ are conditionally independent given $\cF_{t-1}$; moreover,
%
%
\begin{equation} \label{martingale2}
m(\wht\alpha_{\mathrm B}-\alpha) = \sum_{t=1}^{2n-1} \bigl(Z_t^{(1)} +
\cdots+ Z_t^{(m)}\bigr).
\end{equation}

From the martingale representation (\ref{martingale2}) it follows that
$E^*(\wht\alpha_{\mathrm B})=\alpha$. Moreover, under the assumption that
%
%
\begin{equation} \label{finitew}\qquad
\quad\sigma_{\mathrm B}^2 := \sum_{t=1}^n E^* \Biggl[ f_t^2(\bY_t)
\Big/ \prod_{k=1}^{t-1} w_k(\bY_k) \Biggr] E^* \Biggl[ \prod_{k=1}^{t-1}
w_k(\bY_k) \Biggr] -n\alpha^2 < \infty,
\end{equation}
application
of the central limit theorem yields
%
%
\begin{equation} \label{clt}
\sqrt{m} (\wht\alpha_B -\alpha) \Rightarrow N ( 0,
\sigma_{\mathrm B}^2) \qquad\mbox{as } m \rightarrow\infty.
\end{equation}
A consistent
estimate of $\sigma_{\mathrm B}^2$ is given by
%
%
\begin{eqnarray} \label{2.13}
\wht\sigma_{\mathrm B}^2 & := & m^{-1} \sum_{j=1}^m \Biggl\{ \biggl[ \sum_{i\dvtx
a_{n-1}^{(i)}=j} f_n\bigl(\wtd\bY_n^{(i)}\bigr) h_{n-1} \bigl(\bY_{n-1}^{(i)}\bigr)
\biggr] \nonumber\\[-8pt]\\[-8pt]
&&\hphantom{m^{-1} \sum_{j=1}^m \Biggl\{}
{} - \Biggl[1+\sum_{t=1}^{n-1} \sum_{i\dvtx a_{t-1}^{(i)}=j}
\bigl(\#_t^{(i)}- mw_t^{(i)}\bigr) \Biggr] \wht\alpha_{\mathrm B} \Biggr\}^2, \nonumber
\end{eqnarray}
which can be shown to converge to $\sigma_{\mathrm B}^2$ in probability
as $m \rightarrow\infty$ by making use of the martingale representation
(\ref{2.9}) (see \cite{CL08} for details). Del Moral and Jacod \cite
{DJ01} have
derived by a different method a martingale representation similar\vspace*{1pt} to
(\ref{martingale2}) (see \cite{DJ01}, (3.3.7) and (3.3.8)), in which
the term
$Z_{2t-1}^{(i)}$ in~(\ref{list}) corresponds to the $t$th mutation on the
$i$th particle and $Z_{2t}^{(i)}$ the $t$th selection by the $i$th particle.
In~\cite{DJ01}, these two terms are combined into a sum and a central
limit theorem similar to (\ref{clt}) is proved under the assumption of
bounded $f_n$.

Note that in (\ref{clt}) on the asymptotic normality of $\wht\alpha
_{\mathrm B}$
and in the consistency result $\wht\sigma_{\mathrm B}^2 \stackrel
{p}{\rightarrow}
\sigma_{\mathrm B}^2$, the sample size $n$ in the
probability $\alpha= P \{ \bY_n \in\Gamma\}$ is assumed to be
fixed whereas the number $m$ of Monte Carlo samples approaches
$\infty$. The consistent estimate $\wht\sigma_\rB^2$ of
$\sigma_\rB^2$ in (\ref{2.13}) provides an estimate $\wht
\sigma_\rB/\sqrt{m}$ of the standard error (s.e.)($\wht\alpha_\rB
$) of the Monte Carlo
estimate $\widehat\alpha_{\mathrm B}$. Note that the usual estimate
$\sqrt{\wht\alpha_{\rB} (1-\wht\alpha_{\rB})}$ is inconsistent for
$\sqrt m$ s.e.($\wht\alpha_\rB$) because of the dependence among
the $m$ sample paths due to resampling in the SISR procedure as in
\cite{DDJ06,DG05}. The case of $n$ approaching $\infty$ will be
considered in the
next section
in which the representation (\ref{2.10}) will still play
a~pivotal role, but which requires new methods and large deviation
principles rather than central
limit theorems.

Instead of bootstrap resampling, we can use the residual resampling
sche\-me introduced by Baker \cite{Bak85,Bak87} which often leads to
smaller asymptotic variance than that of bootstrap resampling.
We consider here a variant of
this scheme introduced by Crisan, Del Moral and Lyons \cite{CDL99}
that can result in further reduction of the asymptotic
variance. Let $\lfloor\cdot\rfloor$ denote the greatest integer
function and let $m_t$ be the sample size at stage $t$ with $m_1=m$.
We modify the bootstrap resampling step of the SISR procedure as
follows: let $U_t^{(1)},
\ldots, U_t^{(m_t)}$ be independent Bernoulli random variables
satisfying $P \{ U_t^{(i)} = 1 \} = m_t w_t^{(i)} - \lfloor m_t
w_t^{(i)} \rfloor$. For each $1 \leq i \leq m_t$ and $t<n$, make
$\#_t^{(i)}:=\lfloor m_t w_t^{(i)} \rfloor+ U_t^{(i)}$ copies of
$(\wbY_t^{(i)}, a_{t-1}^{(i)},h_{t-1}^{(i)}, w_t^{(i)})$. These
copies constitute an augmented sample $\{ (\bY_t^{(j)},a_t^{(j)},
h_t^{(j)},w_t^{(j)})\dvtx1 \leq j \leq
m_{t+1} \}$, where $m_{t+1} = \sum_{i=1}^{m_t} \#_t^{(i)}$ and
$h_t^{(i)} = h_{t-1}^{(i)}/(m_t w_t^{(i)})$. Estimate $\alpha$ by
\[
\wht\alpha_{\mathrm R} := m_n^{-1} \sum_{i=1}^{m_n} L\bigl(\wbY_n^{(i)}\bigr)
h_{n-1}^{(i)}\bigl(\bY_{n-1}^{(i)}\bigr) {\mathbf1}_{\{ \tilde\bY_n^{(i)} \in
\Gamma\}}.
\]
Define $\varepsilon_k^{(j)}$ by (\ref{2.10}) in which $m$ is replaced
by $m_t$ and define $\cF_{2t-1}$ (or $\cF_{2t}$) by
(\ref{2.11}) in which $m$ is replaced by $m_{s+1}$ (or by $m_{t+1}$).
Moreover, define
\begin{eqnarray*} 
\wtd Z_{2t-1}^{(i)} & = & \bigl[f_t\bigl(\wbY_t^{(i)}\bigr)-f_{t-1}\bigl(\bY_{t-1}^{(i)}\bigr)\bigr]
h_{t-1}\bigl(\bY_{t-1}^{(i)}\bigr) \qquad\mbox{ for } 1 \leq t \leq n, \\
\wtd Z_{2t}^{(i)} & = & \bigl(\#_t^{(i)}-m_t w_t^{(i)}\bigr)\bigl[f_t\bigl(\wbY_t^{(i)}\bigr)
h_t\bigl(\wbY_t^{(i)}\bigr)-
\alpha\bigr] \qquad\mbox{for } 1 \leq t \leq n-1,
\end{eqnarray*}
for $i=1,\ldots,m_t$. Recall that the first generation of particles
consists of $\wtd Y_1^{(1)},\allowbreak\ldots, \wtd Y_1^{(m)}$ and that
$a_t^{(i)}=j$ if the first component of $\bY_t^{(i)}$ is $\wtd Y_1^{(j)}$
for $j=\allowbreak1,\ldots,m$ and $i=1,\ldots,m_{t+1}$.
Analogous to (\ref{2.9}) and (\ref{martingale2}), we
have the martingale representations
%
%
\begin{eqnarray} \label{2.15}
m_n(\wht\alpha_{\mathrm R}-\alpha) &=& \sum_{j=1}^m \bigl(\varepsilon_1^{(j)}+
\cdots+ \varepsilon_{2n-1}^{(j)}\bigr)\nonumber\\[-8pt]\\[-8pt]
&=& \sum_{k=1}^{2n-1} \bigl(\wtd Z_k^{(1)} +
\cdots+ \wtd Z_k^{(m_{\lfloor(k+1)/2 \rfloor})}\bigr).\nonumber
\end{eqnarray}
Analogous to (\ref{2.13}), define
\begin{eqnarray*} 
\wht\sigma_{\mathrm R}^2 & = & m^{-1} \sum_{j=1}^m \Biggl\{ \biggl[ \sum_{i\dvtx
a_{n-1}^{(i)} = j} f_n\bigl(\wtd\bY_n^{(i)}\bigr) h_{n-1}\bigl(\bY_{n-1}^{(i)}\bigr)
\biggr] \\
&&\hphantom{m^{-1} \sum_{j=1}^m \Biggl\{}
{} - \Biggl[ 1+ \sum_{t=1}^{n-1} \sum_{i\dvtx a_{t-1}^{(i)}=j}
\bigl(\#_t^{(i)}-m_t w_t^{(i)}\bigr) \Biggr] \wht\alpha_{\mathrm R} \Biggr\}^2.
\end{eqnarray*}

From (\ref{2.15}) it follows that $E^*[m_n(\wht\alpha_{\mathrm R}-\alpha
)]=0$. Let
\[
\eta_t = E^* \Biggl[ \prod_{k=1}^t
w_k(\bY_k) \Biggr],\qquad h_t^*({\mathbf y}_t)= \eta_t \Big/ \prod_{k=1}^t
w_k({\mathbf y}_k),
\]
and let $\gamma(x)=(x-\lfloor x \rfloor) (1-x+\lfloor x
\rfloor)/x$ for $x > 0$. If (\ref{finitew}) holds,
then analogous to corresponding results for $\wht\alpha_{\mathrm B}$
and $\wht\sigma_{\mathrm B}^2$ in the bootstrap resampling case, we now
have as
$m \rightarrow\infty$,
\begin{eqnarray*}
&\displaystyle \wht\sigma_{\mathrm R}^2 \stackrel{p}{\rightarrow} \sigma_{\mathrm
R}^2,\qquad m_t/m
\stackrel{p}{\rightarrow} 1\qquad
\mbox{for every } t \geq1,&
\\
&\displaystyle \sqrt{m} (\wht\alpha_{\mathrm R}-\alpha) \Rightarrow N(0,\sigma
^2_{\mathrm R}),&
\end{eqnarray*}
where $\sigma_{\mathrm R}^2 < \sigma_{\mathrm B}^2$ and
\begin{eqnarray*} 
\sigma_{\mathrm R}^2 & := & \sum_{t=1}^n E^* \{
[f^2_t(\bY_t)-f_{t-1}^2(\bY_{t-1})] h_{t-1}^*(\bY_{t-1}) \} \\
&&{}
+ \sum_{t=1}^{n-1} E^* \biggl\{ \gamma\biggl( \frac{h_{t-1}^*(\bY
_{t-1})}{h_t^*(\bY_t)}
\biggr) \frac{[f_t(\bY_t) h_t^*(\bY_t) -\alpha]^2}{h_t^*(\bY_t)} \biggr\}.
\end{eqnarray*}
Details are given in \cite{CL08}. Note the
additional variance reduction if residual resampling is
used instead of bootstrap resampling.

\section{Logarithmically efficient SISR for Monte Carlo computation
of small tail probabilities}\label{sec3}

Let
$\xi, \xi_1, \xi_2, \ldots$ be i.i.d. $d$-dimensional random vectors
with a common distribution function $F$ such that $\psi(\theta):=
\log( E e^{\theta' \xi}) < \infty$ for $\| \theta\| < \theta_0$.
Let $S_n=\xi_1+ \cdots+\xi_n$, $\mu_0 = E \xi$, $\Theta= \{
\theta\dvtx
\psi(\theta) < \infty\}$ and let $\Lambda$ be the closure of
$\nabla\psi(\Theta)$ and $\Lambda^o$ be its interior. Assume that
for any $\theta_0 \in\Theta^o$ and $\theta
\in\Theta\setminus\Theta^o$,
\[
\lim_{\rho\uparrow1} (\theta- \theta_0)' \nabla\psi\bigl(\theta_0+
\rho(\theta-\theta_0)\bigr) = \infty.
\]
Then by convex analysis
(see, e.g., \cite{Bro86}, Chapter 3), $\Lambda$ contains the convex
hull of the support of $\{ S_n/n, n \geq1 \}$.
The gradient
vector $\nabla\psi$ is a diffeomorphism from $\Theta^o$ onto
$\Lambda^o$. For given $\mu\in\Lambda^o$ let $\theta_\mu= (\nabla
\psi)^{-1}(\mu)$ and define
the \textit{rate function}
%
%
\begin{equation} \label{3.1}
\phi(\mu) = \sup_{\theta\in\Theta} \{ \theta' \mu-
\psi(\theta) \} = \theta_\mu' \mu-\psi(\theta_\mu).
\end{equation}
We can embed $F$ in an exponential family $\{ F_\theta, \theta\in
\Theta\}$ with
\[
dF_\theta(x) = e^{\theta' x-\psi(\theta)} \,dF(x).
\]
Under certain regularity conditions on $g\dvtx\Lambda\rightarrow\mathbf{R}$,
Chan and Lai \cite{CL00} have developed asymptotic approximations, which
involve both $g$ and $\phi$, to the exceedance probabilities
%
%
\begin{eqnarray} \label{3.2}
p_n & = & P \{ g(S_n/n) \geq b \} \qquad\mbox{with } b > g(\mu_0), \\
\label{3.3}
p_c & = & P \Bigl\{ \max_{n_0 \leq n \leq n_1} n g(S_n/n) \geq c \Bigr\},
\end{eqnarray}
where $n_0 \sim\rho_0 c$ and $n_1 \sim\rho_1 c$ such that
$g(\mu_0) < \rho_1^{-1}$. Making use of these approximations, Chan
and Lai \cite{CL07} have shown that certain mixtures of exponentially
twisted measures are asymptotically optimal for Monte Carlo
evaluation of (\ref{3.2}) or (\ref{3.3}) by importance sampling.
Specifically, for $A = \{ g(S_n/n) \geq b \}$ in the case of
(\ref{3.2}) or $A = \{ \max_{n_0 \leq n \leq n_1} n g(S_n/n) \geq c
\}$ in the case of (\ref{3.3}), an importance sampling measure $Q$
(which may depend on $n$ or $c$) is said to be \textit{asymptotically
optimal} if
%
%
\begin{equation} \label{3.4}
m \var\Biggl( m^{-1} \sum_{i=1}^m L_i {\mathbf1}_{A_i} \Biggr) =
O\bigl(\sqrt{n} p_n^2\bigr) \qquad\mbox{as } n \rightarrow\infty
\end{equation}
in the case of (\ref{3.2}) and if
%
%
\begin{equation} \label{3.5}
m \var\Biggl( m^{-1} \sum_{i=1}^m L_i {\mathbf1}_{A_i} \Biggr) =
O(p_c^2) \qquad\mbox{as } c \rightarrow\infty
\end{equation}
in the case of (\ref{3.3}), where $(L_1,{\mathbf1}_{A_1}), \ldots,
(L_m, {\mathbf1}_{A_m})$ are $m$ independent realizations of
($L:=dP/dQ$, ${\mathbf1}_A$). For the case of (\ref{3.3}), since $E_Q(L
{\mathbf1}_A) = P(A)=p_c$, $E_Q(L^2 {\mathbf1}_A) \geq p_c^2$ by the
Cauchy--Schwarz inequality and, therefore, $Q$ is an asymptotically
optimal importance sampling measure if $E_Q(L^2 {\mathbf1}_A) =
O(p_c^2)$, which leads to the definition (\ref{3.5}) of asymptotic
optimality for the Monte Carlo estimates. Chan and Lai \cite{CL07} have
also shown that $\sqrt{n} p_n^2$ is an asymptotically minimal order
of magnitude for $E_Q(L^2 {\mathbf1}_A)$ in the case of (\ref{3.2}).
They have also extended this theory to Markov random walks $S_n$
whose increments~$\xi_i$ have distributions $F(\cdot|X_i, X_{i-1})$
depending on a Markov chain~$\{ X_t \}$.

The asymptotically optimal mixtures of exponentially twisted
measures $\int P_{\theta_\mu} \omega(\mu) \,d \mu$ in \cite{CL07}
involve normalizing constants $\beta_n$ (or $\beta_c$) that
may be difficult to compute. Moreover, it may even be difficult to
sample from the twisted measure $P_{\theta_\mu}$, especially in
multidimensional and Markovian settings. In this section we show
that by choosing the resampling weights suitably, the SISR estimates
$\widehat\alpha_{\rB}$ can still attain
%
%
\begin{equation} \label{3.6}
m \var(\widehat\alpha_{\rB}) = p_n^2 e^{o(n)} \qquad\mbox{ as } m
\rightarrow\infty\mbox{ and } n \rightarrow\infty
\end{equation}
for Monte Carlo estimation of $p_n$ and
%
%
\begin{equation} \label{3.7}
m \var(\widehat\alpha_{\rB}) =p_c^2 e^{o(c)} \qquad\mbox{as } m
\rightarrow\infty\mbox{ and } c \rightarrow\infty
\end{equation}
for Monte Carlo estimation of (\ref{3.3}). Moreover, (\ref{3.6}) and
(\ref{3.7}) still hold with~$\wht\alpha_\rB$ replaced by $\wht
\alpha_{\mathrm R}$. The properties (\ref{3.6}) and (\ref{3.7}) are
called \textit{logarithmic efficiency}; the variance of the Monte Carlo
estimate differs from the asymptotically optimal value by a factor
of $e^{o(n)}$ (or $e^{o(c)}$) noting that $-n^{-1} \log p_n$ and
$-c^{-1} \log p_c$ converge to positive limits. To begin with,
suppose the asymptotically optimal importance sampling measure $Q$
has conditional densities $q_t(\cdot|\bY_{t-1})$ with respect to
$\nu$. To achieve log efficiency,\vadjust{\goodbreak} the resampling functions $w_t$ can
be chosen to satisfy approximately
%
%
\begin{equation} \label{wtby}
w_t({\mathbf y}_t) \propto q_t(y_t|{\mathbf y}_{t-1})/ \widetilde
q_t(y_t|{\mathbf y}_{t-1})
\end{equation}
as illustrated by the following example, after which a heuristic explanation
for (\ref{wtby}) will be given.
\begin{exa}\label{exam1}
Suppose $\xi_1, \xi_2, \ldots$ are i.i.d. random
variables ($d=1$) and $g(x)=x$ in (\ref{3.2}), so that $\alpha= p_n = P
\{ S_n/n \geq b \}$, where $b> E \xi_1$ and $2 \theta_b \in\Theta$.
Consider the SISR procedure with $\wQ=P$ (and, therefore, $E^*=E$) and
resampling weights
%
%
\begin{equation} \label{3.8}
w_t(\bY_t) = e^{\theta_b \xi_t-\psi(\theta_b)}.
\end{equation}
Then $L=1$ and hence, by (\ref{2.5}),
%
%
\begin{equation} \label{3.9}
f_t(\bY_t) = P\{ S_n/n \geq b |\bY_t \} = P \{S_n - S_t \geq n b -
S_t | S_t\}.
\end{equation}
Therefore, standard Markov's inequality involving moment generating
functions yields
%
%
\begin{equation} \label{3.10}
f_t(\bY_t) \leq e^{-\theta_b (nb - S_t) + (n-t) \psi(\theta_b)}
= e^{\theta_b S_t-t \psi(\theta_b) - n
\phi(b)}.
\end{equation}
By (\ref{2.10}) and the martingale decomposition (\ref{2.9}),
%
%
\begin{eqnarray} \label{3.11}
E(\widehat\alpha_\rB- \alpha)^2
&\leq&
m^{-1} \sum_{t=1}^n E \bigl\{\bigl[f_t\bigl(\wbY_t^{(1)}\bigr)-f_{t-1}\bigl(\bY_{t-1}^{(1)}\bigr)\bigr]^2
h_{t-1}^2\bigl(\bY_{t-1}^{(1)}\bigr) \bigr\} \nonumber\\[-8pt]\\[-8pt]
&&{} + m^{-1} \sum_{t=1}^{n-1}
E \bigl[\bigl(
\#_t^{(1)} - mw_t^{(1)}\bigr)^2 f_t^2\bigl(\wbY_t^{(1)}\bigr)
h_t^2\bigl(\wbY_t^{(1)}\bigr)\bigr], \nonumber
\end{eqnarray}
in which the superscript $^{(1)}$ can be replaced by $^{(i)}$ since the
expectations
are the same for all $i$.
The derivation of (\ref{3.11}) uses the independence of
$[f_t(\wbY_t^{(i)})- f_{t-1}(\bY_{t-1}^{(i)})] h_t(\bY_{t-1}^{(i)})$
for $1 \leq i \leq m$ when conditioned on $\cF_{2t-2}$ and the
pairwise negative correlations of $(\#_t^{(i)}-mw_t^{(i)}) f_t(\wbY_t^{(i)})
h_t(\wbY_t^{(i)})$ for $i=1,\ldots,m$ when conditioned on
$\cF_{2t-1}$. By (\ref{2.4a}), (\ref{3.8}) and (\ref{3.10}),
%
%
\begin{eqnarray} \label{3.12}
&& E \bigl\{ \bigl[ f_t\bigl(\wbY_t^{(1)}\bigr)-
f_{t-1}\bigl(\bY_{t-1}^{(1)}\bigr)\bigr]^2 h_{t-1}^2\bigl(\bY_{t-1}^{(1)}\bigr) \bigr\} \nonumber\\
&&\qquad
= E \bigl\{ \bar w_1^2 \cdots\bar w_{t-1}^2
\bigl[f_t\bigl(\wbY_t^{(1)}\bigr)-f_{t-1}\bigl(\bY_{t-1}^{(1)}\bigr)\bigr]^2/ e^{2 \theta_b
S_{t-1}^{(1)}- 2(t-1) \psi(\theta_b)} \bigr\} \\
&&\qquad \leq\biggl(
1+\frac{E(e^{\theta_b \xi_1- \psi(\theta_b)}-1)^2}{m} \biggr)^{t-1}
e^{-2n \phi(b)} E\bigl(e^{2 \theta_b \xi_t-2 \psi(\theta_b)}\bigr).
\nonumber
\end{eqnarray}
To see the inequality in (\ref{3.12}), condition on $\cF_{2t-1}$. Since
$E[f_t(\wtd\bY_t^{(1)})|\cF_{2t-1}]=f_{t-1}(\bY_{t-1}^{(1)})$,
it follows from (\ref{3.10}) that
\begin{eqnarray*}
&& E \bigl\{ \bigl[f_t\bigl(\wbY_t^{(1)}\bigr)-f_{t-1}\bigl(\bY_{t-1}^{(1)}\bigr)\bigr]^2/e^{2 \theta
_b S_{t-1}^{(1)}
-2(t-1) \psi(\theta_b)} | \cF_{2t-1} \bigr\}\\
&&\qquad \leq E\bigl[f_t^2\bigl(\wbY_t^{(1)}\bigr)/e^{2
\theta_b S_{t-1}^{(1)}-2(t-1)\psi(\theta_b)}|\cF_{2t-1}\bigr] \leq
e^{-2n \phi(b)} E\bigl(e^{2 \theta_b \xi_t-2 \psi(\theta_b)}\bigr).
\end{eqnarray*}
Moreover, $\bar w_1^2, \ldots, \bar w_{t-1}^2$ are i.i.d. random variables
with mean
%
%
\begin{equation} \label{resample}\quad
\quad E \Biggl[ m^{-1} \sum_{i=1}^m \bigl(e^{\theta_b \xi_1^{(i)}-\psi(\theta_b)}-1\bigr)+1
\Biggr]^2
= 1+m^{-1} E\bigl(e^{\theta_b \xi_1-\psi(\theta_b)}-1\bigr)^2
\end{equation}
and their product $\bar w_1^2 \cdots\bar w_{t-1}^2$ in the second term of
(\ref{3.12}) is $\cF_{2t-1}$-measurable. This yields the inequality
in (\ref{3.12}).

Since the conditional distribution of $\#_t^{(i)}$ given $\cF_{2t-1}$
is Binomial$(m,w_t^{(i)})$, $E[(\#_t^{(i)}-mw_t^{(i)})^2|\cF_{2t-1}]
\!\leq\!mw_t^{(i)}$.
By (\ref{2.4a}), (\ref{3.8}) and (\ref{3.10}), $f_t(\wbY_t^{(i)})
h_t(\wbY_t^{(i)})\!{\leq}\allowbreak\bar w_1 \cdots\bar w_t e^{-n \phi(b)}$. Since $\sum_{i=1}^m
w_t^{(i)}\,{=}\,1$,
it then follows by conditioning on~$\mathcal{F}_{2t-1}$ that
\begin{eqnarray*}
&& E \bigl\{ \bigl(\#_t^{(1)}-mw_t^{(1)}\bigr)^2 f_t^2\bigl(\wbY_t^{(1)}\bigr)
h_t^2\bigl(\wbY_t^{(1)}\bigr) \bigr\} \\
&&\qquad = m^{-1} \sum_{i=1}^m
E \bigl\{ \bigl(\#_t^{(i)}-mw_t^{(i)}\bigr)^2 f_t^2\bigl(\wbY_t^{(i)}\bigr)
h_t^2\bigl(\wbY_t^{(i)}\bigr) \bigr\} \\
&&\qquad \leq E \Biggl\{ \Biggl( \sum_{i=1}^m w_t^{(i)} \Biggr)\bigl(\bar w_1 \cdots\bar w_t
e^{-n \phi(b)}\bigr)^2 \Biggr\} = e^{-2n \phi(b)} E(\bar w_1^2 \cdots\bar
w_t^2),
\end{eqnarray*}
which can be combined with (\ref{resample}) to yield
%
%
\begin{equation} \label{3.13}
\qquad E\bigl[\bigl(\#_t^{(1)}-mw_t^{(1)}\bigr)^2 f_t^2\bigl(\wbY_t^{(1)}\bigr)
h_t^2\bigl(\wbY_t^{(1)}\bigr)\bigr]
= O \biggl( \biggl(1+ \frac{K}{m} \biggr)^t e^{-2n \phi(b)} \biggr),
\end{equation}
where $K = E(e^{\theta_b \xi_1- \psi(\theta_b)}-1)^2$. By
(\ref{3.11}), (\ref{3.12}) and (\ref{3.13}),
\[
\liminf_{n \rightarrow\infty} - \frac{1}{n} \log
[m \var(\widehat\alpha_\rB)] \geq2 \phi(b) - \frac{K}{m}
\]
for any fixed $m$. Since $p_n/[n^{-1/2} e^{-n \phi(b)}]$ is bounded
away from 0 and $\infty$
(see~\cite{CL07}, page 451), (\ref{3.6}) holds.
\end{exa}

\subsection{A heuristic principle for efficient SISR procedures}\label{sec3.1}

The asymptotically optimal importance\vspace*{1pt} sampling measure for $p_n = P \{ S_n/n
\geq b \}$ is $Q$ under which $\xi_1, \xi_2, \ldots$ are i.i.d. with
density function $e^{\theta_b \xi- \psi(\theta_b)}$ with respect
to $P$ (see \cite{CL07}). Since we have used $\wtd Q=P$
in Example \ref{exam1}, (\ref{3.8}) actually follows the prescription
(\ref{wtby})
to choose resampling weights that can achieve an effect similar
to asymptotically optimal importance sampling. We now give a heuristic
principle underlying this prescription. The SISR procedure uses the
importance weights $p_t^{(i)}/\wtd q_t^{(i)}$ (for the change of
measures from $P$ to $\wtd Q$) and resampling weights $w_t^{(i)}$,
$1 \leq i \leq m$, for the $m$ simulated trajectories at stage $t$.
The resampling\vspace*{1pt} step at stage $t$ basically converts $(\wbY_t^{(i)},
p_t^{(i)}/\widetilde q_t^{(i)},w_t^{(i)})$ to $(\bY_t^{(i)},p_t^{(i)}/
(\widetilde q_t^{(i)} w_t^{(i)}),1)$, and, therefore, the prescription
(\ref{wtby}) for choosing resampling weights (satisfying $\widetilde q_t^{(i)}
w_t^{(i)} = q_t^{(i)}$) is intended to yield
the desired importance weights $p_t^{(i)}/q_t^{(i)}$.
To transform this heuristic principle into a rigorous proof of
logarithmic efficiency, one needs to be able to bound the second
moments of the importance weights and resampling weights. This
explains the requirement $2 \theta_b \in\Theta$ in Example \ref{exam1}.

Example \ref{exam1} indicates the key role played by the martingale
decomposition~(\ref{2.9}) and large deviation bounds for $P(\Gamma_n
| \mathbf{Y}_k)$, $1 \leq k < n$, in the derivation of asymptotically
efficient resampling weights. To generalize the basic ideas to the
more general tail probability (\ref{3.2}) with nonlinear $g$, we provide
large deviation bounds in Lemma \ref{lemma1}, whose proof is given in
the \hyperref[app]{Appendix}, for
%
%
\begin{equation} \label{3.14}
P \bigl\{g\bigl((x + S_{n, k}) / n\bigr) \geq b\bigr\},
\end{equation}
where $S_{n, k} = S_n - S_k$; note that (\ref{3.14}) is equal to
$P\{g (S_n / n) \geq b | S_k = x\}$. The special case $k=0$ and $x =
0$ has been analyzed by Chan and Lai (see Theorem 2 of \cite{CL00})
under certain
regularity conditions that yield precise saddlepoint approximations.
The probability (\ref{3.14}) is more complicated than this special
case because it involves
additional parameters $x$ and $k$, but we only need large deviation
bounds rather than saddlepoint approximations for logarithmic
efficiency. Let $\mu_\theta= \nabla\psi(\theta)$ and define
%
%
\begin{eqnarray} \label{I}
I & = & \inf\{ \phi(\mu)\dvtx
g(\mu) \geq b \}, \\
\label{Htheta}
M & = & \{ \theta\dvtx\phi(\mu_\theta) \leq I \}.
\end{eqnarray}

\begin{lem}\label{lemma1}
Let $b > g(\mu_0)$. Then as $n \rightarrow\infty$,
%
%
\begin{equation} \label{3.16}
P \bigl\{ g\bigl((x+S_{n,k})/n\bigr) \geq b \bigr\} \leq e^{-n I+o(n)} \int_M
e^{\theta' x-k \psi(\theta)} \,d \theta,
\end{equation}
where the $o(n)$ term is
uniform in $x$ and $k$.
\end{lem}

The proof of (\ref{3.16}) in the \hyperref[app]{Appendix} uses a change-of-measure
argument that involves the measure $Q$ for which
\[
(dQ/dP)({\mathbf Y}_n)
= \int_M e^{\theta' S_n-n \psi(\theta)}
\,d \theta/\operatorname{vol}(M).
\]
The bound (\ref{3.16})
is used in conjunction with the inequality $\int_M e^{\theta' x-k
\psi(\theta)}
\,d \theta\leq\operatorname{vol}(M) \exp\{ k \max_{\theta\in M} [\theta'
x/k-\psi(\theta)] \}$
to prove the following theorem.
\begin{theorem}\label{theo1} Letting $b > g(\mu_0)$, assume:
\begin{longlist}[(C1)]
\item[(C1)]
$g$ is twice continuously differentiable and
$\nabla g \not= 0$ on
$N : = \{\mu\in\Lambda^o \dvtx g (\mu)
= b\}$.
\item[(C2)]
$Ee^{2 \kappa\| \xi_1 \|} < \infty$,
where $\kappa=
{\sup_{\theta\in M}} \| \theta\|$
and $M$ is defined in (\ref{Htheta}).
\end{longlist}
Let $\wht\theta_0=0$ and define for $1 \leq t \leq n$,
%
%
\begin{eqnarray} \label{3.17}
\wht\theta_t & = & \mathop{\arg\max}_{\theta\in M} \{ \theta'
S_t/t-\psi(\theta) \}, \nonumber\\[-8pt]\\[-8pt]
w_t(\bY_t) & = & \exp\{ \wht\theta'_t S_t
-t \psi(\wht\theta_t)-[\wht\theta'_{t-1} S_{t-1} -(t-1) \psi(\wht
\theta_{t-1})] \}. \nonumber
\end{eqnarray}
With $\wQ=P$ and the resampling weights thus defined, the SISR
estimates~$\wht\alpha_\rB$ and $\wht\alpha_\rR$ are
logarithmically efficient, that is, (\ref{3.6}) holds for $\wht
\alpha_\rB$ and also with~$\wht\alpha_\rR$ in place of $\wht
\alpha_\rB$ if $m \rightarrow\infty$ and $n \rightarrow\infty$.
\end{theorem}

Besides (\ref{3.16}), the proof of Theorem \ref{theo1} also uses the
bounds in the
following lemma. These bounds enable us to bound $E(\bar
w_{t-1}^2|\mathcal{F}_{2(t-1)-2})$
in the proof of Theorem \ref{theo1}.
\begin{lem}\label{lemma2} With the same notation and
assumptions in Theorem \ref{theo1}, there exist nonrandom constants
$\varepsilon_t$ and $K > 0$ such that
%
%
\begin{eqnarray} \label{6--1}
&\displaystyle \lim_{t \rightarrow\infty} \varepsilon_t = 0,\qquad
E[ w_t(\bY_t)|S_{t-1} ] \leq e^{\varepsilon_t} \quad\mbox{and
}&\nonumber\\[-8pt]\\[-8pt]
&\displaystyle E[w_t^2(\bY_t)|S_{t-1}] \leq K \qquad\mbox{for all } t
\geq1.&\nonumber
\end{eqnarray}
\end{lem}
\begin{pf}
Let $\eta= {\sup_{\theta\in M}} |\psi(\theta)|$. Then
%
%
\begin{eqnarray} \label{A.16}
\wht\theta_t' S_t-t \psi(\wht\theta_t) & = & [\wht\theta_t' S_{t-1}-(t-1)
\psi(\wht\theta_t)]+[\wht\theta_t' \xi_t-\psi(\wht\theta_t)]
\nonumber\\[-8pt]\\[-8pt]
& \leq&
[\wht\theta_{t-1}' S_{t-1}-(t-1) \psi(\wht\theta_{t-1})]+[\wht
\theta_t'
\xi_t-\psi(\wht\theta_t)] \nonumber
\end{eqnarray}
and, therefore, it follows from (\ref{3.17}) that $w_t(\bY_t) \leq
e^{\kappa\|
\xi_t \|+\eta}$. Hence, by (C2),
%
%
\begin{equation} \label{A.17}\qquad
\quad E\bigl[w_t(\bY_t) {\mathbf1}_{\{ \| \xi_t \| > \zeta\}}|
S_{t-1}\bigr] \leq E\bigl[e^{\kappa\| \xi_1 \|+\eta} {\mathbf1}_{\{ \| \xi_1 \|
> \zeta
\}}\bigr] \rightarrow0 \qquad\mbox{as } \zeta\rightarrow\infty.
\end{equation}
It will be shown that for any fixed $\zeta> 0$,
%
%
\begin{equation} \label{A.18}
\gamma_{t,\zeta} := {\operatorname{ess} \sup}\| \wht
\theta_t -\wht\theta_{t-1} \| {\mathbf1}_{\{ \| \xi_t \| \leq\zeta\}}
\rightarrow0 \qquad\mbox{as } t \rightarrow\infty.
\end{equation}
Let $\wtd\eta= {\sup_{\theta\in M}} \| \nabla\psi(\theta) \|$.
Combining (\ref{A.18}) with (\ref{3.17}) and (\ref{A.16}) yields
%
%
\begin{eqnarray} \label{A.19}
E\bigl[w_t(\bY_t) {\mathbf1}_{\{ \| \xi_t \| \leq\zeta\}}|S_{t-1}\bigr] &\leq&
E\bigl[e^{\hat\theta_t' \xi_t-\psi(\hat\theta_t)} {\mathbf1}_{\{ \| \xi
_t \|
\leq\zeta\}}|S_{t-1}\bigr] \nonumber\\
&\leq& e^{\gamma_{t,\zeta} (\zeta+\tilde\eta)}
E\bigl[e^{\hat\theta_{t-1}' \xi_t-\psi(\hat\theta
_{t-1})}|S_{t-1}\bigr]\\
&=&1+o(1) \nonumber
\end{eqnarray}
as $t \rightarrow\infty$. Moreover, by (C2) and (\ref{A.18}), as
$\zeta\rightarrow\infty$,
%
%
\begin{eqnarray} \label{wt2}
\quad E\bigl[w_t^2(\bY_t) {\mathbf1}_{\{ \| \xi_t \| > \zeta\}}|S_{t-1}\bigr] & \leq&
E\bigl[e^{2 \kappa\| \xi_1 \|+2 \eta} {\mathbf1}_{\{ \| \xi_1 \| > \zeta\}}\bigr]
\rightarrow0 \nonumber\\
\quad E\bigl[w_t^2(\bY_t) {\mathbf1}_{\{ \| \xi_t \| \leq\zeta\}}| S_{t-1}\bigr]
&\leq&
e^{2 \gamma_{t,\zeta} (\zeta+\tilde\eta)}
E\bigl[e^{2 \hat\theta_{t-1}' \xi_t-2 \psi(\hat\theta_{t-1})}|S_{t-1}\bigr]
\\
\quad & \leq& \sup_{\theta\in M} e^{\psi(2 \theta)-2 \psi(\theta
)}+o(1), \nonumber
\end{eqnarray}
and (\ref{6--1}) follows from (\ref{A.17}), (\ref{A.19}) and\vadjust{\goodbreak} (\ref{wt2}).

To prove (\ref{A.18}), let $f_{x,t}(\theta)=\theta' x-t \psi(\theta
)$ and let
$\theta_{x,t}$ be the unique maximizer of $f_{x,t}(\theta)$
over $M$. Let $\lambda_{\min}(\cdot)$ denote the smallest eigenvalue
of a~symmetric matrix. Since $\nabla^2 \psi(\theta)$ is continuous and
positive definite for
all $\theta\in M$, and since $M$ is compact and $\lambda_{\min}$ is
a continuous
function of the entries of $\nabla^2 \psi(\theta)$, $\inf_{\theta
\in M}
\lambda_{\min}(\nabla^2 \psi(\theta)) \geq2 \beta$ for some
$\beta> 0$.
Therefore, by Taylor's theorem,
$f_{x,t-1}(\theta) \leq f_{x,t-1}(\theta_{x,t-1}) -\beta t \|
\theta_{x,t-1} - \theta\|^2$ for all $\theta\in M$. It then follows that
for $\| y-x \| \leq\zeta$,
\begin{eqnarray*}
f_{y,t}(\theta_{x,t-1}) &\leq& f_{y,t}(\theta_{y,t})
= f_{x,t-1}(\theta_{y,t})+ \theta_{y,t}'(y-x)-\psi(\theta_{y,t})
\\
&\leq& f_{x,t-1}(\theta_{x,t-1})-\beta t \|
\theta_{x,t-1}
-\theta_{y,t} \|^2 + \theta_{y,t}'(y-x)-\psi(\theta_{y,t}) \\
&\leq& f_{y,t}(\theta_{x,t-1})-\beta t \| \theta_{x,t-1} -
\theta_{y,t} \|^2 +(\zeta+\wtd\eta) \| \theta_{x,t-1}-\theta
_{y,t} \|
\end{eqnarray*}
and, therefore, $\| \theta_{x,t-1} - \theta_{y,t} \| \leq
(\zeta+\wtd\eta)/(\beta t)$. Hence, (\ref{A.18}) holds by setting
$x=S_{t-1}$ and $y=S_t$.
\end{pf}
\begin{pf*}{Proof of Theorem \ref{theo1}} To simplify the notation, we
will suppress
the superscript $^{(1)}$ in $\wht\theta_{t-1}^{(1)}$ below. By
(\ref{2.4a}) and (\ref{3.17}),
%
%
\begin{equation}\label{3.18}
h_{t-1}\bigl(\wbY_{t-1}^{(1)}\bigr) = \Biggl( \prod_{k=1}^{t-1} \bar w_k \Biggr)
\exp\bigl[-\wht\theta_{t-1}' \wtd S_{t-1}^{(1)}+(t-1) \psi(\wht
\theta_{t-1})\bigr].
\end{equation}
Making use of $E[f(\wbY_t^{(1)})|\cF_{2t-2}]\,{=}\,f_{t-1}(\bY_{t-1}^{(1)})$,
$E(\sup_{\theta\in M} e^{2 \theta' \xi_t-2 \psi(\theta)})\,{<}\,
\infty$~and the independence of
$\bar w_1^2 \cdots\bar w_{t-1}^2$ and $\xi_t$, we obtain from
Lemma \ref{lemma1} and~(\ref{3.18}) that
%
%
\begin{eqnarray} \label{3.19}\quad
&& E \bigl\{ \bigl[f_t\bigl(\wbY_t^{(1)}\bigr)-f_{t-1}\bigl(\bY_{t-1}^{(1)}\bigr)\bigr]^2
h_{t-1}^2\bigl(\bY_{t-1}^{(1)}\bigr) \bigr\} \nonumber\\
&&\qquad \leq E \bigl\{ \bar w_1^2 \cdots
\bar w_{t-1}^2 f_t^2\bigl(\wbY_t^{(1)}\bigr)/\exp\bigl[2 \wht\theta_{t-1}'
S_{t-1}^{(1)}-2(t-1) \psi(\wht\theta_{t-1})\bigr] \bigr\} \\
&&\qquad \leq
e^{-2n
I+o(n)} E(\bar w_1^2 \cdots\bar w_{t-1}^2). \nonumber
\end{eqnarray}
By (\ref{2.4a}) and Lemma \ref{lemma2},
\begin{eqnarray*}
E\bigl(\bar w_{t-1}^2|\cF_{2(t-1)-2}\bigr) &=& \Biggl( m^{-1} \sum_{i=1}^m
E\bigl[w_{t-1}\bigl(\wbY_{t-1}^{(i)}\bigr) |S_{t-2}^{(i)}\bigr] \Biggr)^2 \\
&&{} +m^{-2}
\sum_{i=1}^m \var\bigl[w_{t-1}\bigl(\wbY_{t-1}^{(i)}\bigr)| S_{t-2}^{(i)}\bigr] \\
&\leq&
(1+Km^{-1}) e^{2 \varepsilon_{t-1}}
\end{eqnarray*}
and proceeding inductively yields
%
%
\begin{equation} \label{3.20}\qquad
\quad E(\bar w_1^2 \cdots\bar w_{t-1}^2) \leq(1+K m^{-1})^{t-1} \exp\Biggl(
\sum_{k=1}^{t-1} 2 \varepsilon_k \Biggr) \leq e^{K (t-1)/m+o(n)}.
\end{equation}
Similarly, under bootstrap or residual resampling,
%
%
\begin{eqnarray} \label{3.21}
&& E\bigl[\bigl(\#_t^{(1)}-mw_t^{(1)}\bigr)^2
f_t^2\bigl(\wbY_t^{(1)}\bigr)
h_t^2\bigl(\wbY_t^{(1)}\bigr)\bigr] \nonumber\\
&&\qquad = m^{-1} \sum_{i=1}^m
E\bigl[\bigl(\#_t^{(i)}-mw_t^{(i)}\bigr)^2 f_t^2\bigl(\wbY_t^{(i)}\bigr)
h_t^2\bigl(\wbY_t^{(i)}\bigr)\bigr] \\
&&\qquad \leq e^{-2nI+o(n)} E(\bar w_1^2
\cdots\bar w_t^2). \nonumber
\end{eqnarray}
By (C1), $p_n = e^{-nI+o(n)}$ (see \cite{CL00}, Theorem 2) and hence,
it follows from~(\ref{3.11}) and (\ref{3.19})--(\ref{3.21}) that both $\wht
\alpha_{\mathrm R}$ and $\wht\alpha_{\mathrm B}$ are logarithmically efficient.
\end{pf*}

The heuristic principle described in the paragraph following
Example \ref{exam1} can also be used to construct logarithmically
efficient SISR
procedures for Monte Carlo evaluation of (\ref{3.3}) as illustrated in
the following example.
\begin{exa}\label{exam2}
Let $T_c = \inf\{ n\dvtx S_n \geq c \}$. Consider the
estimation of $p_c = P \{ T_c \leq n_1 \}$ [i.e., with $d=1$ and $g(x)=x$]
when $\mu_0 <0$ and $n_1 \sim ac$ for
some $a > 1/ \psi'(\theta_*)$, where $\theta_*$ is the unique
positive root of $\psi(\theta_*)=0$. We shall assume $2 \theta_*
\in\Theta$ and use the importance measure $\wQ=P$ and resampling weights
\[
w_t(\bY_t) =
\cases{e^{\theta_* \xi_t}, &\quad if $t \leq T_c$,
\cr
1, &\quad if $n_1 > t > T_c$.}
\]
Let $\eta(\bY_{T_c \wedge n_1}) = e^{\theta_*(S_{T_c \wedge
n_1}-c)}$. Since
$\eta(\bY_{T_c \wedge n_1}) \geq{\mathbf1}_{\{ \max_{n \leq n_1} S_n
\geq c \}}$, it
follows that
%
%
\begin{equation} \label{3.14a}
\qquad f_t(\bY_t) = P \Bigl\{ \max_{n \leq n_1} S_n \geq c \big| \bY_t \Bigr\}
\leq E[\eta(\bY_{T_c \wedge n_1})|\bY_t] = e^{\theta_* (S_{T_c
\wedge
t}-c)}.
\end{equation}
Making use of (\ref{3.14a}) in place of (\ref{3.10}), we obtain that,
analogous to (\ref{3.12}),
%
%
\begin{eqnarray} \label{3.15a}
&&E \bigl\{ \bigl[ f_t\bigl(\wbY_t^{(1)}\bigr)-f_{t-1}\bigl(\bY_{t-1}^{(1)}\bigr)\bigr]^2 h_{t-1}^2
\bigl(\bY_{t-1}^{(1)}\bigr) \bigr\} \nonumber\\[-8pt]\\[-8pt]
&&\qquad\leq\biggl( 1 +\frac{K_*}{m} \biggr)^{t-1} e^{-2
\theta_* c} E(e^{2 \theta_* \xi_t}),\nonumber
\end{eqnarray}
where $K_* = E(e^{\theta_* \xi_1}-1)^2$ and that, analogous to (\ref{3.13}),
%
%
\begin{equation} \label{3.16a}\qquad
\quad E\bigl[\bigl(\#_t^{(1)}-mw_t^{(1)}\bigr)^2 f_t^2\bigl(\wbY_t^{(1)}\bigr) h_t^2\bigl(\wbY_t^{(1)}\bigr)\bigr]
=O \biggl( \biggl( 1 + \frac{K_*}{m} \biggr)^{t-1} e^{-2 \theta_* c}
\biggr).
\end{equation}
Hence, by (\ref{3.11}) (with $n_1$ in place of $n$), (\ref{3.15a})
and (\ref{3.16a}),
\[
m \var(\widehat\alpha_\rB) = O \bigl( n_1 \exp[ (n_1 K_*/m) - 2
\theta_* c ] \bigr).
\]
Since $n_1=O(c)$ and $p_c/ e^{-\theta_* c}$ is bounded away from 0 and
$\infty$,
as shown in~\cite{Sie75}, (\ref{3.7}) also holds.
\end{exa}

In Theorem \ref{theo2}, we provide the resampling weights for logarithmically
efficient simulation of (\ref{3.3}), for which the counterparts of
(\ref{I}) and
(\ref{Htheta}) are also provided. The basic idea is to use the resampling
weights (\ref{3.17}) up to the stopping time
%
%
\begin{equation} \label{3.31}
T_c = \inf\{ n \geq n_0\dvtx n g(S_n/n) \geq c \} \wedge n_1.
\end{equation}

\begin{theorem}\label{theo2}
$\!\!$Let $g(\mu_0)\,{<}\,a^{-1}$, $n_0\,{=}\,\delta
c+O(1)$ and $n_1\,{=}\,ac+O(1)$
as $c\,{\rightarrow}\,\infty$ for some $a > \delta>
0$. Let $I = \inf\{ \phi(\mu)\dvtx g(\mu) \geq\delta^{-1} \}$
and\vspace*{1pt}
$M = \{ \theta\dvtx\phi(\mu_\theta) \leq I
\}$. Let $\wQ= P$ and assume that \textup{(C1)--(C2)} hold for all $a^{-1}
\leq b \leq\delta^{-1}$ and that

\begin{longlist}[(C3)]
\item[(C3)]
$r:= \sup_{\mu\dvtx g(\mu) \geq a^{-1}} \min\{ g(\mu),
\delta^{-1} \}/\phi(\mu)
< \infty$.
\end{longlist}
Let $\wht\theta_0=0$ and define for $1 \leq t
\leq n_1-1$, $\wht\theta_t = \arg\max_{\theta\in M} [\theta'
S_t/t-\psi(\theta)]$ and
%
%
\begin{equation} \label{thetahat}\qquad
w_t(\bY_t) = \cases{ e^{ \hat\theta_t' S_t-t \psi(\wht\theta_t)
- [\hat
\theta_{t-1}' S_{t-1} -
(t-1) \psi(\wht\theta_{t-1})]}, &\quad if $t \leq T_c$,\cr
1, &\quad if $n_1 > t > T_c$.}
\end{equation}
Then (\ref{3.7}) holds for $\wht\alpha_{\mathrm B}$ and with
$\wht\alpha_{\mathrm B}$ replaced by $\wht\alpha_{\mathrm R}$ if $m
\rightarrow\infty$ and $c \rightarrow\infty$.
\end{theorem}
\begin{pf}
Let $u=(t-1) \wedge T_c^{(1)}$. By (\ref{2.4a})
and (\ref{thetahat}),
%
%
\begin{equation} \label{hseq}
h_{t-1}\bigl(\wbY_{t-1}^{(1)}\bigr) = \Biggl( \prod_{k=1}^{t-1} \bar w_k \Biggr) \exp\bigl[
-\bigl(\wht\theta_u^{(1)}\bigr)' \wtd S_u^{(1)}+u \psi\bigl(\wht\theta_u^{(1)}\bigr)\bigr].
\end{equation}
Let $I_b = \inf\{ \phi(\mu)\dvtx g(\mu) \geq b \}$. By Lemma \ref{lemma1},
%
%
\begin{eqnarray} \label{ftseq}\qquad\qquad
f_t\bigl(\wbY_t^{(1)}\bigr) &=& P \bigl\{ T_c^{(1)} \leq n_1 | \wbY_t^{(1)} \bigr\}
\nonumber\\[-8pt]\\[-8pt]
&\leq&
\cases{\displaystyle \sum_{n=t+1}^{n_1} e^{-n I_{c/n}+o(n)} \int_M
e^{\theta' \tilde S_t^{(1)}-t \psi(\theta)} \,d \theta, &\quad if $t
< T_c^{(1)}$, \vspace*{2pt}\cr
1, &\quad if $t \geq T_c^{(1)}$.}
\nonumber
\end{eqnarray}
Note that
\[
\inf_{a^{-1} \leq b \leq\delta^{-1}} b^{-1} I_b =
\min\biggl\{ \inf_{\mu\dvtx a^{-1}
\leq g(\mu) \leq\delta^{-1}} \frac{\phi(\mu)}{g(\mu)}, \inf
_{\mu\dvtx
g(\mu) > \delta^{-1}} \frac{\phi(\mu)}{\delta^{-1}} \biggr\} = r^{-1}
\]
by (C3). Hence, by (\ref{hseq}) and (\ref{ftseq}),
%
%
\begin{equation} \label{Fseq}
E \bigl\{ \bigl[f_t\bigl(\wbY_t^{(1)}\bigr)-f_{t-1}\bigl(\bY_{t-1}^{(1)}\bigr)\bigr]^2 h_{t-1}^2\bigl(\bY
_{t-1}^{(1)}\bigr) \bigr\}
\leq e^{-2c/r+o(c)} E(\bar w_1^2 \cdots\bar w_{t-1}^2).\hspace*{-35pt}
\end{equation}
Similarly, it can be shown that under either bootstrap or residual resampling,
%
%
\begin{equation} \label{hexseq}\quad
E\bigl[\bigl(\#_t^{(1)}-mw_t^{(1)}\bigr)^2 f_t^2\bigl(\wbY_t^{(1)}\bigr) h_t^2\bigl(\wbY_t^{(1)}\bigr)\bigr]
\leq e^{-2c/r+o(c)} E(\bar w_1^2 \cdots\bar w_{t-1}^2).
\end{equation}
By (C1) and Theorem 2 of \cite{CL00}, $p_c=e^{-c/r+o(c)}$ and hence,
it follows from (\ref{3.20}), (\ref{Fseq})
and (\ref{hexseq}) that both $\wht\alpha_{\mathrm R}$ and $\wht\alpha
_{\mathrm B}$
are logarithmically efficient.
\end{pf}

\subsection{Markovian extensions}\label{sec3.2}

Let $\{ (X_t,S_t)\dvtx t=0,1,\ldots, \}$ be a Markov additive process on
$\cX\times\mathbf{R}^d$ with transition kernel
\begin{eqnarray*}
P(x,A \times B) :\!&=& P \{ (X_1,S_1) \in A \times(B+s) |
(X_0,S_0)=(x,s) \} \\
&=& P \{ (X_1,S_1) \in A \times B|
(X_0,S_0) = (x,0) \}.
\end{eqnarray*}
Let $\{ X_n \}$ be aperiodic and irreducible with respect to some
maximal irreducibility measure $\varphi$ and assume that the
transition kernel satisfies the minorization condition
%
%
\begin{equation} \label{m1}
P(x,A \times B) \geq h(x,B) \nu(A)
\end{equation}
for any measurable set $A \subset\cX$, Borel set $B \subset
\mathbf{R}^d$ and $s \in{\mathbf R}^d$ for some probability measure $\nu
$ and
measure $h(x,\cdot)$ that is positive for all $x$ belonging to a
$\varphi$-positive set. Ney and Nummelin \cite{NN87} developed a theory
to analyze large deviations properties of $S_n$ under (\ref{m1}) or
when its variant $P(x,A \times B) \geq h(x) \nu(A \times B)$ holds.
Let $\tau$ be the first regeneration time and assume that
$\Omega:= \{ (\theta,\zeta)\dvtx E_\nu e^{\theta' S_\tau-\tau\zeta} <
\infty\}$ is an open neighborhood of 0.
Then for all $\theta\in\Theta:= \{ \theta\dvtx(\theta,\zeta) \in
\Omega$ for some $\zeta\}$, the kernel
%
%
\begin{equation} \label{whtPtheta}
\wht P_\theta(x,A) := \int
e^{\theta' s} P(x,A \times ds)
\end{equation}
has a unique maximum eigenvalue
$e^{\psi(\theta)}$, for which $\zeta=\psi(\theta)$ is the unique
solution of the equation
$E_\nu e^{\theta' S_\tau-\tau\zeta}=1$, with corresponding right
eigenfunctions $r(\cdot; \theta)$ and left eigenmeasures $\ell_\nu
(\theta, \cdot)$ defined by
%
%
\begin{eqnarray} \label{rxtheta}
r(x;\theta) & = & E_x e^{\theta' S_\tau-\tau
\psi(\theta)}, \nonumber\\
\ell_x(\theta;A) & = & E_x \Biggl( \sum_{n=0}^{\tau-1}
e^{\theta' S_n-n \psi(\theta)} {\mathbf1}_{\{ X_n \in A \}} \Biggr),
\\
\ell_\nu(\theta; A) & = & \int\ell_x(\theta; A) \,d \nu(x).
\nonumber
\end{eqnarray}
Let $\pi$ denote the stationary distribution of $\{
X_n \}$ and let
%
%
\begin{equation} \label{thetamu}
\theta_\mu= (\nabla\psi)^{-1}(\mu).
\end{equation}

To begin with,
consider the special case $d=1$ and $g(x)=x$ for which the importance
sampling measure with transition kernel
%
%
\begin{equation} \label{Ptheta}
P_\theta(x,dy \times ds) := e^{\theta' s-\psi(\theta)}
\{ r(y;\theta)/r(x;\theta) \} P(x, dy \times ds)
\end{equation}
has been shown to be
logarithmically efficient by Dupuis and Wang \cite{DW05} and
asymptotically optimal by Chan and Lai \cite{CL07} for simulating the tail
probability $P_{x_0} \{ S_n/n \geq b \}$ when $\theta$ is chosen to be
$\theta_b$ in (\ref{Ptheta}). We shall show that by
using SISR with $\wtd Q=P$ and resampling
weights $w_t(\bY_t)=e^{\theta_b \xi_t-\psi(\theta_b)}$, we can
avoid computation
of the eigenfunctions. To bring out the essence of the method,
we first assume instead of the minorization
condition (\ref{m1}) the stronger uniform recurrence condition
%
%
\begin{equation} \label{ur}
a_0 \nu(A \times B) \leq P(x,A \times B) \leq a_1 \nu(A \times B)
\end{equation}
for some $0 < a_0 < a_1$ and probability measure $\nu$ and for all $x
\in\cX$,
measurable sets $A \subset\cX$ and Borel sets $B \subset{\mathbf R}$.
At the end of this section, we
show how this assumption can be removed. Note that $\bY_t$ consists of
$(X_i,\xi_i)$, $i \leq t$,
in the Markov case.

\begin{exa}\label{exam3}
Let $b> E_\pi\xi_1$ and assume that $\theta_b \in
\Theta$ and $E_\nu(e^{2 \theta_b \xi_1-2 \psi(\theta_b)})$ $<
\infty$.
We now extend Example \ref{exam1} to Markov additive processes by
showing that the
choice $\wtd Q=P$ and
%
%
\begin{equation} \label{weights}
w_t(\bY_t) = e^{\theta_b \xi_t-\psi(\theta_b)}
\end{equation}
results in logarithmically efficient simulation of
$P_{x_0} \{ S_n/n \geq b \}$. The dependence of the weights $w_t^{(i)}$
and $w_t^{(j)}$ for $i \neq j$, created from a combination of the Markovian
structure of the underlying process and bootstrap resampling, requires
a more
delicate peeling and induction argument than that in Example \ref
{exam1}. By considering
$\xi_t-\psi(\theta_b)/\theta_b$ instead of $\xi_t$, we may assume
without loss
of generality that $\psi(\theta_b)=0$.

Let $\kappa= \sup_{x \in\cX} r(x;\theta_b)/\inf_{x \in\cX}
r(x;\theta_b)$
and let $E_\theta$ be expectation with respect to $P_\theta$.
Then by (\ref{2.5}) and (\ref{Ptheta}),
\begin{eqnarray*}
f_t(\bY_t) & = & P_{x_0} \{ S_n/n \geq b | \bY_t \} = P
\{ S_n-S_t \geq nb-S_t|X_t,S_t \} \\
& = & r(X_t;\theta_b) E_{\theta_b} \bigl[e^{-\theta_b(S_n-S_t)}
\mathbf{1}_{\{ S_n-S_t \geq nb-S_t \}}/
r(X_n; \theta_b)|X_t, S_t\bigr] \\
& \leq& \kappa e^{-\theta_b(nb-S_t)}.
\end{eqnarray*}
We shall show that
%
%
\begin{equation} \label{moment}\quad
E(\bar w_1^2 \cdots\bar w_t^2) = e^{o(t)} \qquad\mbox{as } m \rightarrow
\infty
\mbox{ uniformly over } 1 \leq t \leq n-1.
\end{equation}
Then logarithmic efficiency of bootstrap resampling follows from
(\ref{3.11})--(\ref{3.13}). We first show that for any $k<t$ and $i
\neq j$,
%
%
\begin{eqnarray} \label{peel}
&& E \bigl\{ \bar w_k^2 \bigl(E_{X_k^{(i)}} e^{\theta_b S_{t-k}}\bigr)\bigl(
E_{X_k^{(j)}} e^{\theta_b S_{t-k}}\bigr)|
\cF_{2k-2} \bigr\} \nonumber\\[-8pt]\\[-8pt]
&&\qquad \leq m^{-2} \sum_{u \neq v} \bigl(E_{X_{k-1}^{(u)}} e^{\theta_b S_{t-k+1}}\bigr)
\bigl(E_{X_{k-1}^{(v)}} e^{\theta_b S_{t-k+1}}\bigr) + m^{-1} \beta, \nonumber
\end{eqnarray}
where $\beta=\sup_{h \geq0, x \in\cX} E_x \{ e^{2 \theta_b \xi
_1} (E_{X_1} e^{\theta_b
S_h})^2 \}$, which is finite by (\ref{ur}). Note that
$\bar w_k$ is measurable with respect to $\cF_{2k-1}$
and that under bootstrap resampling, $X_k^{(i)}$ and $X_k^{(j)}$ are
independent conditioned
on $\cF_{2k-1}$. Moreover, since $X_k^{(1)} = \wtd X_k^{(\ell)}$ with
probability $w_k^{(\ell)} =
w_k(\wtd\bY^{(\ell)}_k)/\sum_{j=1}^m w_k(\wtd\bY_k^{(j)})$,
\[
E \bigl\{ \bar w_k \bigl(E_{X_k^{(1)}}e^{\theta_b S_{t-k}}\bigr)|\cF_{2k-1} \bigr\} =
\bar w_k \sum_{u=1}^m
w_k^{(u)} E_{\tilde X_k^{(u)}} e^{\theta_b S_{t-k}},
\]
which is equal to $m^{-1} \sum_{u=1}^m
e^{\theta_b \tilde\xi_k^{(u)}} E_{\tilde X_k^{(u)}} e^{\theta_b
S_{t-k}}$ in view of
(\ref{weights}) and that $\psi(\theta_b)=0$. Hence,
%
%
\begin{eqnarray} \label{step1}
&& E \bigl\{ \bar w_k^2 \bigl(E_{X_k^{(i)}} e^{\theta_b S_{t-k}}\bigr)\bigl(
E_{X_k^{(j)}} e^{\theta_b S_{t-k}}\bigr)
| \cF_{2k-1} \bigr\} \nonumber\\
&&\qquad = \Biggl(m^{-1} \sum_{u=1}^m e^{\theta_b \tilde\xi_k^{(u)}}
E_{\tilde X_k^{(u)}} e^{\theta_b S_{t-k}} \Biggr)^2
\nonumber\\[-8pt]\\[-8pt]
&&\qquad = m^{-2} \sum_{u \neq v} \bigl(e^{\theta_b \tilde\xi_k^{(u)}}
E_{\tilde
X_k^{(u)}} e^{\theta_b S_{t-k}}\bigr)\bigl( e^{\theta_b \tilde\xi_k^{(v)}}
E_{\tilde X_k^{(v)}} e^{\theta_b S_{t-k}}\bigr) \nonumber\\
&&\qquad\quad{} + m^{-2} \sum_{u=1}^m e^{2 \theta_b \tilde\xi_k^{(u)}}
\bigl(E_{\tilde X_k^{(u)}}
e^{\theta_b S_{t-k}}\bigr)^2. \nonumber
\end{eqnarray}
Since $(\wtd\xi_k^{(u)},\wtd X_k^{(u)})$ and $(\wtd\xi_k^{(v)},\wtd
X_k^{(v)})$ are
independent conditioned on $\cF_{2k-2}$ for $u \neq v$ and
$E[e^{\theta_b \tilde
\xi_k^{(i)}}(E_{\tilde X_k^{(i)}} e^{\theta_b S_{t-k}})| \cF_{2k-2}]
= E_{X_{k-1}^{(i)}}
e^{\theta_b S_{t-k+1}}$, (\ref{peel}) follows from (\ref{step1}).

We shall show using (\ref{peel}) and induction, that
%
%
\begin{equation} \label{bound}\qquad
E(\bar w_1^2 \cdots\bar w_k^2) \leq\gamma^2 (1+m^{-1} \beta)^k
\quad\mbox{where }
\gamma= \sup_{x \in\cX, h \geq0} E_x e^{\theta_b S_h} (\mbox{$\geq$}1).
\end{equation}
For $k=1$,
\[
E \bar w_1^2 = m^{-2} \sum_{i \neq j} E_{x_0} e^{\theta_b \xi_1^{(i)}}
E_{x_0} e^{\theta_b \xi_1^{(j)}} +m^{-2} \sum_{i=1}^m E_{x_0} e^{2
\theta_b
\xi_1^{(i)}} \leq\gamma^2 +m^{-1} \beta
\]
and indeed (\ref{bound}) holds. If (\ref{bound}) holds for all $k <
t$, then
by repeated application of (\ref{peel}), starting from $k=t$, we obtain
\begin{eqnarray*}
E(\bar w_1^2 \cdots\bar w_t^2) & \leq& (E_{x_0} e^{\theta_b S_t})^2+m^{-1}
\beta\sum_{k=0}^{t-1} E(\bar w_1^2 \cdots\bar w_k^2) \\[-2pt]
& \leq& \gamma^2 \Biggl\{ 1+m^{-1} \beta
\sum_{k=0}^{t-1} (1+m^{-1} \beta)^k \Biggr\} = \gamma^2 (1+m^{-1} \beta
)^t
\end{eqnarray*}
and (\ref{bound}) indeed holds for $k=t$. Hence, (\ref{moment}) is true
and logarithmic efficiency is attained.
\end{exa}

The peeling argument used to derive (\ref{peel}) and (\ref{bound})
can also be used to extend Theorems \ref{theo1} and \ref{theo2}, which
hold for general
$g$, to the following.\vadjust{\goodbreak}
\begin{theorem}\label{thm3}
\textup{(a)}
Let $M$, $\wht\theta_t$ and $w_t(\bY_t)$ be the same as in
Theorem~\ref{theo1}. Then
Theorem \ref{theo1} still holds when the i.i.d. assumption on $\xi_t$ is
replaced by
the uniform recurrence condition (\ref{ur}) on the Markov
additive process
$(X_t,S_t=\xi_1+\cdots+\xi_t)$ and assumption \textup{(C2)} is
generalized to
%
%
\begin{equation} \label{expo}
\int_{{\mathbf R}^d} e^{2 \kappa\| \xi\|} \nu(\cX, d \xi) < \infty
\qquad\mbox{where }
\kappa= {\sup_{\theta\in M}} \| \theta\|.
\end{equation}

\textup{(b)} Let $M$, $\wht\theta_t$ and $w_t(\bY_t)$ be the same as in
Theorem \ref{theo2}. Then
Theorem~\ref{theo2} still holds when the i.i.d. assumption on $\xi_t$ is
replaced by
the uniform recurrence condition (\ref{ur}) and assumption
\textup{(C2)} is generalized to
(\ref{expo}).\vspace*{-3pt}
\end{theorem}

Note that $\wtd Q=P$ in Theorem \ref{thm3}.
We next show how the uniform recurrence assumption (\ref{ur}) can be
removed, extending the
preceding results on the logarithmic efficiency of suitably chosen
SISR procedures to more general Markov additive
processes such that for some $\theta\in\Theta$, $0 < \beta< 1$,
function $u\dvtx\cX
\rightarrow[1,\infty)$ and measurable set $C$:

\begin{longlist}[(U1)]
\item[(U1)]
$\sup_{x \in C} u(x) < \infty$, $\int
_\mathcal{X} u(x) \,d \nu(x) < \infty$,
$\sup_{x \in C} \ell_x(\theta;C) < \infty$, $\int
_\mathcal{X} \ell_x(\theta$; $C) \,d \nu(x) < \infty$,\vspace*{1pt}

\item[(U2)] $E_x \{ e^{\theta' \xi_1-\psi(\theta)} u(X_1) \} \leq(1-\beta
) u(x)$
for $x \notin C$,\vspace*{1pt}

\item[(U3)] $a:=\sup_{x \in C} E_x \{ e^{\theta' \xi_1-\psi(\theta)}
u(X_1) \}
< \infty$,

\item[(U4)] $K_1:= \sup_{x \in\cX} E_x \{ e^{2 \theta' \xi_1-2 \psi
(\theta)} u^2(X_1)
/u^2(x) \} < \infty$.
\end{longlist}
We illustrate in Section \ref{sec4}, Example \ref{exam5}, how (U1)--(U4)
can be checked in a~concrete
example. Condition (U1) [in which $\ell_x$ is defined in (\ref
{rxtheta})] holds when $C$ is bounded and
$\nu$ has support on a compact set. Conditions (U2)--(U4) are often
called ``drift conditions''
(see \cite{CL07}). Although the arguments are essentially
modifications of the
peeling idea in Example \ref{exam3} by making use of (U1)--(U4), they are
considerably more complicated than those in the uniformly
recurrent case. We, therefore, only consider the univariate linear case
[$d=1$, $g(y)=y$]
in the following theorem to
indicate the basic ideas without getting into the details
of these modifications, such as replacing for general $g$ the
$\theta_b$ in (\ref{uweights}) by sequential estimates $\wht\theta
_t$, as in (\ref{3.17})
and (\ref{thetahat}).\vspace*{-3pt}

\begin{theorem}\label{theo4}
Let $b > E_\pi\xi_1$ and
assume that \textup{(U1)--(U4)} hold for $\theta=\theta_b$. Let $\wtd
Q=P$ and
%
%
\begin{equation} \label{uweights}
w_t(\bY_t) = e^{\theta_b \xi_t-\psi(\theta_b)} u(X_t)/u(X_{t-1}).
\end{equation}
Then (\ref{3.6}) holds with $p_n = P_{x_0} \{ S_n/n \geq b \}$, for
$\wht\alpha_{\mathrm B}$ or
$\wht\alpha_{\mathrm R}$, as $n \rightarrow\infty$ and $m \rightarrow
\infty$.
\end{theorem}
\begin{pf}
By considering $\xi_t-\psi(\theta_b)/\theta_b$ instead of $\xi_t$,
we assume without loss of generality that $\psi(\theta_b)=0$. By
(\ref{2.4a})
and (\ref{uweights}),
%
%
\begin{equation} \label{hu}
h_{t-1}\bigl(\wbY_{t-1}^{(1)}\bigr) = \Biggl( \prod_{k=1}^{t-1} \bar w_k \Biggr)
e^{-\theta_b \tilde S_{t-1}^{(1)}} u(x_0)/u\bigl(\wtd
X_{t-1}^{(1)}\bigr).\vadjust{\goodbreak}
\end{equation}
It will be shown in the \hyperref[app]{Appendix} that
%
%
\begin{equation} \label{ku}
K_2:= \sup_{x \in\cX, h \geq0} E_x \{ e^{\theta_b S_h} u(X_h)/u(x)
\} < \infty.
\end{equation}
Note that
%
%
\begin{eqnarray} \label{fu}
f_t(\bY_t) & = & E_{x_0}\bigl({\mathbf1}_{\{ S_n/n \geq b \}}|\bY_t\bigr) \leq
e^{-\theta_b
nb} E_{x_0} (e^{\theta_b S_n}|\bY_t) \nonumber\\[-8pt]\\[-8pt]
& = & e^{\theta_b (S_t-nb)} E_{X_t}(e^{\theta_b S_{n-t}}) \leq K_2
e^{\theta_b
(S_t-nb)} u(X_t). \nonumber
\end{eqnarray}
Since $E_{x_0} [f_t(\wbY_t^{(1)})|\cF_{2t-2}]\,{=}\,f_{t-1}(\bY
_{t-1}^{(1)})$, it follows from
(\ref{hu}), (\ref{fu}) and~(U3) that
%
%
\begin{eqnarray} \label{sequ}\quad
&& E_{x_0} \bigl\{ \bigl[f_t\bigl(\wbY_t^{(1)}\bigr)-f_{t-1}\bigl(\bY_{t-1}^{(1)}\bigr)\bigr]^2
h_{t-1}^2\bigl(\bY_{t-1}^{(1)}\bigr) \bigr\}
\nonumber\\
&&\qquad \leq K_2^2 e^{-2n \theta_b b} E_{x_0} \bigl\{ (\bar w_1 \cdots\bar w_{t-1})^2
e^{2 \theta_b \tilde\xi_t^{(1)}} u^2(x_0) u^2\bigl(\wtd
X_t^{(1)}\bigr)/u^2\bigl(X_{t-1}^{(1)}\bigr) \bigr\} \\
&&\qquad \leq \beta e^{-2n \theta_b b} E_{x_0} (\bar w_1^2 \cdots\bar
w_{t-1}^2), \nonumber
\end{eqnarray}
where $\beta=K_1 K_2^2 u^2(x_0)$.

By (\ref{hu}) and (\ref{fu}), under either bootstrap or residual resampling,
%
%
\begin{eqnarray} \label{resampleu}
&& E_{x_0} \bigl[\bigl(\#_t^{(1)}-mw_t^{(1)}\bigr)^2 f_t^2\bigl(\wbY_t^{(1)}\bigr)h_t^2\bigl(\wbY
_t^{(1)}\bigr)\bigr] \nonumber\\
&&\qquad = m^{-1} \sum_{i=1}^m E_{x_0} \bigl[\bigl(\#_t^{(i)}-mw_t^{(i)}\bigr)^2
f_t^2\bigl(\wbY_t^{(i)}\bigr)
h_t^2\bigl(\wbY_t^{(i)}\bigr)\bigr] \\
&&\qquad \leq K_2^2 E_{x_0} (\bar w_1^2 \cdots\bar w_t^2) e^{-2n \theta_b b}
u^2(x_0). \nonumber
\end{eqnarray}
In view of (\ref{3.11}), it now remains to show (\ref{moment}). It follows
from the proof of (\ref{peel}) that for any $k<t$ and $i \neq j$,
\begin{eqnarray*} 
&& E_{x_0} \biggl\{ \bar w_k^2 \biggl(\frac{E_{X_k^{(i)}} [e^{\theta_b S_{t-k}}
u(X_{t-k})]}{
u(X_k^{(i)})} \biggr) \biggl(\frac{E_{X_k^{(j)}} [e^{\theta_b S_{t-k}} u(X_{t-k})]}{
u(X_k^{(j)})} \biggr) \Big| \cF_{2k-2} \biggr\} \\
&&\qquad \leq m^{-2} \sum_{v \neq w}
\biggl(\frac{E_{X_{k-1}^{(v)}} [e^{\theta_b S_{t-k+1}} u(X_{t-k+1})]}{
u(X_{k-1}^{(v)})} \biggr)\!\biggl(\frac{E_{X_{k-1}^{(w)}} [e^{\theta_b S_{t-k+1}}
u(X_{t-k+1})]}{
u(X_{k-1}^{(w)})} \biggr) \\
&&\qquad\quad{} + m^{-1} \beta.
\end{eqnarray*}
An argument similar to that in (\ref{peel}) and (\ref{bound}) can be
used to show
that
\[
E_{x_0} (\bar w_1^2 \cdots\bar w_k^2) \leq K_2^2 (1+m^{-1} \beta)^k.
\]
Hence, (\ref{moment}) again holds and (\ref{3.6}) follows from (\ref
{sequ}) and
(\ref{resampleu}).
\end{pf}

\subsection{Implementation, estimation of standard errors and
discussion}\label{sec3.3}

As explained in the first paragraph of Section \ref{sec3.1}, at every
stage $t$,
the SISR procedure
carries out importance sampling sequentially within each simulated
trajectory but performs
resampling across the $m$ trajectories. Since the computation\vadjust{\goodbreak} time for
resampling increases
with $m$, it is more efficient to divide the $m$ trajectories into $r$
subgroups of size $k$
so that $m=kr$ and resampling is performed within each subgroup of $k$
trajectories,
independently of the other subgroups. This method also has the
advantage of providing a direct
estimate of the standard error of the Monte Carlo estimate $\bar\alpha
:= r^{-1}
\sum_{i=1}^r \wht
\alpha_i$, where $\wht\alpha_i$ denotes the SISR estimate of $\alpha
$ (using either bootstrap
or residual resampling) based on the $i$th subgroup of simulated
trajectories. Specifically,
we can estimate the standard error of $\bar\alpha$ by $\wht\sigma
/\sqrt{r}$, where
%
%
\begin{equation} \label{sub}
\wht\sigma^2 = (r-1)^{-1} \sum_{i=1}^r (\wht\alpha_i - \bar\alpha)^2.
\end{equation}
In Section \ref{sec2} we considered the case of fixed $n$ as $m
\rightarrow
\infty$ and provided
estimates of the standard errors of the asymptotically normal $\wht
\alpha_{\mathrm B}$ and
$\wht\alpha_{\mathrm R}$. The validity of these estimates is unclear for
the case $n \rightarrow
\infty$ and $m \rightarrow\infty$ as considered in this section that
involves large deviations
theory instead of central limit theorems. By choosing $m=kr$ with $k
\rightarrow\infty$ and
$r \rightarrow\infty$ in~(\ref{sub}), we still have a consistent
estimate $\wht
\sigma/\sqrt{r}$ of the standard error in the large deviations
setting with $n \rightarrow
\infty$.

The resampling weights in Theorems \ref{theo1} and \ref{theo2} have
closed-form expressions
in terms of the
cumulant generating function $\psi(\theta)$ in the i.i.d. case or the
logarithm
$\psi(\theta)$ of the largest eigenvalue of the kernel (\ref{whtPtheta})
in the Markov case. When $\psi(\theta)$
does not have an explicit formula, we can use numerical approximations
and thereby approximate
the logarithmically efficient resampling weights, as will be
illustrated in Example \ref{exam5}.
This is, therefore, much more flexible than logarithmically efficient
importance sampling
which involves sampling from the efficient importance measure that
involves both
the eigenvalue and
corresponding eigenfunction in the Markov case (see
\cite{BNS90,CL07,Col02,DW05,SB90}).
Note that
approximating the eigenvalue and eigenfunction usually does not result
in an importance
(probability) measure and, therefore, requires an additional task of
computing the normalizing
constants.

The basic ideas in Examples \ref{exam1} and \ref{exam2} and Sections
\ref{sec3.1} and \ref{sec3.2} can be
extended to more general
rare events of the form $\{ \bX_T \in\Gamma\}$ and more general
stochastic sequences $\bX_t$
and stopping times $T$. To evaluate $P \{ \bX_T \in\Gamma\}$ by
Monte Carlo, it would be
ideal to sample from the importance measure $Q$ for which
%
%
\begin{equation} \label{3.60}
\frac{dQ}{dP}(\bX_t) = P \{ \bX_T \in\Gamma| \bX_t \}/P \{ \bX_T
\in\Gamma\} \qquad\mbox{for }
t \leq T,
\end{equation}
because the corresponding Monte Carlo estimate
of $P \{ \bX_T \in\Gamma\}$ would have variance 0 (see \cite{DW05},
page 2).
This is clearly not feasible because the
right-hand side of (\ref{3.60}) involves
the conditional probabilities $P \{ \bX_T \in\Gamma| \bX_t \}$ and
its expectation
$P \{ \bX_T \in\Gamma\}$ which is an unknown quantity to be
determined. On the other hand, SISR
enables one to ignore the normalizing factor $P \{ \bX_T \in\Gamma\}$
and to use tractable approximations to
$P \{ \bX_T \in\Gamma| \bX_t \}$, as in Example \ref{exam1}, in coming up
with a logarithmically efficient
Monte Carlo estimate of $P \{ \bX_T \in\Gamma\}$.

\section{Illustrative examples}\label{sec4}

We use the following two examples to
illustrate Theorems \ref{theo1} and \ref{theo4}.
\begin{exa}\label{exam4}
Let $X_1, X_2, \ldots$ be i.i.d. random variables
with $EX_1=0$. Let $\xi_i=(X_i,X_i^2)$ and $S_n = \xi_1 + \cdots+
\xi_n$. Define
$g(y,v)=y/\sqrt{v}$ for $y \in{\mathbf R}$ and $v>0$ and note that
$g(S_n/n)$ is the
self-normalized sum of the $X_i$'s. There is extensive literature on the
large deviation probability $p_n = P \{ g(S_n/n) \geq b \}$ (see \cite{DLS09}).
Consider the case $b=1/\sqrt{2}$ and $X_i$ having the density
function
\[
f(x) = \frac{1}{2 \sqrt{2 \pi}}
\bigl(e^{-(x-1)^2/2}+e^{-(x+1)^2/2}\bigr),\qquad
x \in{\mathbf R},
\]
with respect to the Lebesgue measure. Thus, $X_i$ is a mixture of
$N(1,1)$ and
$N(-1,1)$. In this case, $\Theta= \{ (\theta_1,\theta_2)\dvtx\theta_2
< 1/2 \}$,
$\Lambda= \{ (y,v)\dvtx v \geq y^2 \}$ and
\[
\log(Ee^{\theta_1 X_1+\theta_2 X_1^2}) = \log\biggl( \frac{1}{2} \biggr) +
\frac{1}{2}
-\frac{\theta_1^2+1}{2-4 \theta_2}
+ \log\biggl(\frac{e^{\theta_1/(1-2 \theta_2)} + e^{-\theta_1/(1-2
\theta_2)
}}{\sqrt{1-2 \theta_2}} \biggr)
\]
for $\theta\in\Theta$.
The infimum of the rate function over the one-dimensional
manifold $N = \{ (y,v): y = \sqrt{v/2} \}$ is $I=0.324$ and is
attained at $(y,v)=(1,2)$. Then $M = \{ \theta=(\theta_1,
\theta_2)\dvtx
\phi(y_\theta,v_\theta) \leq I \}$ [see (\ref{Htheta}) and Theorem \ref{theo1}].
We implement SISR with bootstrap resampling as described in Section \ref
{sec3.3},
with $m=10$,000 particles, divided into 100 groups each having 100
particles. The results, in
the form of mean${}\pm{}$standard error and for $n=15, 20$ and 25, are
summarized
in Table \ref{table1}, which also compares them to corresponding results obtained
%
%
\begin{table}[b]
\tablewidth=280pt
\caption{Monte Carlo estimates of $P \{ g(S_n/n) \geq1/\sqrt{2}
\}$}\label{table1}
\begin{tabular*}{\tablewidth}{@{\extracolsep{\fill}}lccc@{}}
\hline
& \multicolumn{3}{c@{}}{$\bolds{n}$}\\[-4pt]
& \multicolumn{3}{c@{}}{\hrulefill}\\
& \textbf{15} & \textbf{20} & \textbf{25} \\
\hline
SISR & $(1.10 \pm0.07) \times10^{-3}$ & $(1.9\pm0.2)\times10^{-4}$
& $(4.0\pm0.7) \times10^{-5}$
\\
Direct & $(0.9 \pm 0.3) \times10^{-3}$ & $(1\pm1) \times10^{-4}$ &
0 \\
\hline
\end{tabular*}
\end{table}
by direct Monte Carlo
with $m=10$,000 in (\ref{2.1}) and (\ref{2.2}). Table \ref{table1} shows 18-fold
variance reduction by
using SISR when $n=15$, 25-fold variance reduction when $n=20$ and that
direct Monte Carlo
fails when $n=25$.
\end{exa}
\begin{exa}\label{exam5}
Let $\zeta_1, \zeta_2, \ldots, \gamma_1, \gamma_2, \ldots$ be
i.i.d. standard normal random variables and let
%
%
\begin{equation} \label{nonlinear}
X_{n+1} = \lambda(X_n) + \zeta_{n+1}, \qquad\xi_n = X_n + \gamma_n,
\end{equation}
where $\lambda(x)$ is a monotone increasing, piecewise linear function
given by
\[
\lambda(x) = x {\mathbf1}_{\{ |x| \leq1 \}}+ \biggl( \frac{x+1}{2} \biggr)
\mathbf{1}_{\{ x>1 \}}
+ \biggl( \frac{x-1}{2} \biggr) {\mathbf1}_{\{ x<-1 \}}.
\]
Let $\theta> 0$. We now show that
(U1)--(U4) hold for $u(x) = e^{2.1 \theta x^+}$ and $C=(-\infty,\rho
]$, where
$\rho\geq1$ is chosen large enough so that (U2) holds, as shown below.
Since $(a+b)^+ \leq a+b^+$ for
$a > 0$ and since $e^{2.05 \theta x} \leq e^{-0.05 \theta x} u(x)$,
it follows that for $x > \rho$,
\begin{eqnarray*}
E_x \bigl\{ e^{\theta\xi_1-\psi(\theta)} u(X_1) \bigr\} & = & E \bigl\{ e^{\theta
x +
\theta\gamma_1 - \psi(\theta) + 2.1 \theta(({x+1})/{2}+\zeta
_1)^+} \bigr\} \\
& \leq& u(x) e^{-0.05 \theta x} E \bigl\{ e^{\theta\gamma_1 - \psi(\theta)
+1.05 \theta+2.1 \theta\zeta_1^+} \bigr\}
\end{eqnarray*}
and, therefore, (U2) holds if $\rho$ is large enough. It is easy to
check that\vspace*{1pt} (U3) holds. Note
that $\sup_{x \in(-\infty,1]} E_x[e^{2 \theta\xi_1-2 \psi(\theta)}
u^2(X_1)] < \infty$ and that for $x>1$,
\begin{eqnarray*}
&& E_x\bigl[e^{2 \theta\xi_1-2 \psi(\theta)} u^2(X_1)\bigr]/u^2(x) \\
&&\qquad = E \bigl[e^{2 \theta x+2 \theta\gamma_1 - 2 \psi(\theta)+4.2 \theta
(({x+1})/{2}+
\zeta_1)^+}\bigr]/e^{4.2 \theta x^+} \\
&&\qquad \leq (e^{-0.1 \theta x} \wedge e^{2 \theta x}) E \bigl[e^{2 \theta\gamma
_1 -2 \psi(\theta)+ 2.1 \theta+4.2 \theta
\zeta_1^+}\bigr] \rightarrow0 \qquad\mbox{as } x \rightarrow\infty
\end{eqnarray*}
and, therefore, (U4) holds. Since $\lim_{x \rightarrow-\infty}
E_x(e^{\theta\xi_1-\psi(\theta)})=0$, it follows
that $\lim_{x \rightarrow-\infty} \ell_x(\theta;C)=0$; moreover,
$u(x)=1$ for all $x \leq0$
and hence, (U1) also holds.

We compute $P_0 \{ S_n/n \geq2.5 \}$ for SISR using
resampling, with $m=10$,000 particles divided into 100 groups, each
having 100 particles, and with resampling weights (\ref{uweights}) for
which the following
procedure is used to provide a numerical approximation for $\theta
_{2.5}$. First note that by
(\ref{nonlinear}),
%
%
\begin{equation} \label{Exetheta}
E_x e^{\theta\xi_1} = e^{\theta^2/2} E_x e^{\theta X_1}.
\end{equation}
The procedure involves a finite-state Markov chain approximation to
(\ref{nonlinear}) with
states $x_i$ and transition probabilities $p_{ij}$ ($1 \leq i, j \leq
1\mbox{,}000$) given by
\[
x_i = \frac{i}{100}-2.505,\qquad p_{ij} = e^{-(x_j-\lambda(x_i))^2/2} \Big/
\sum_{k=1}^{1\mbox{,}000} e^{-(x_k-\lambda(x_i))^2/2}.
\]
For given $\theta$, it approximates $\psi(\theta)$ by $\theta^2/2 +
\wtd\psi(\theta)$, where
$e^{\tilde\psi(\theta)}$ is the largest eigenvalue of the matrix
$(e^{\theta x_j} p_{ij})_{1 \leq i,j \leq1\mbox{,}000}$, in view of (\ref
{whtPtheta}) and
(\ref{Exetheta}). Since $\psi'(\theta_{2.5})=2.5$ by (\ref
{thetamu}), it uses Brent's method
\cite{PFTV92} that involves bracketing followed by\vspace*{1pt} efficient search to
find the positive
root $\wtd\theta_{2.5}$ of the equation $\wtd\psi(\theta)+\theta^2/2
= 2.5 \theta$, noting that $\wtd\psi(0)=0$. The root $\wtd\theta
_{2.5}=0.273$ is then
used as an approximation to $\theta_{2.5}$ in (\ref{uweights}). Table
\ref{table2} gives the results,
%
%
\begin{table}
\caption{Monte Carlo estimates of $P_0 \{ S_n/n \geq2.5 \}$}\label{table2}
\begin{tabular*}{\tablewidth}{@{\extracolsep{\fill}}ld{1.3}ccc@{}}
\hline
& & \multicolumn{3}{c@{}}{$\bolds{n}$} \\[-4pt]
& & \multicolumn{3}{c@{}}{\hrulefill}\\
& \multicolumn{1}{c}{$\bolds\theta$} & \textbf{15} & \textbf{20} & \textbf{25}
\\
\hline
SISR & 0.1 & $(9.68\pm1.37) \times10^{-4}$ & $(2.81\pm0.57) \times
10^{-4}$ &
$(4.70\pm1.22) \times10^{-5}$ \\
& 0.2 & $(9.65\pm0.75) \times10^{-4}$ & $(2.45\pm0.24) \times10^{-4}$ &
$(6.70\pm0.64) \times10^{-5}$ \\
& 0.273 & $(8.31\pm0.48) \times10^{-4}$ & $(2.42\pm0.19) \times
10^{-4}$ &
$(6.33\pm0.44) \times10^{-5}$ \\
& 0.3 & $(9.11\pm0.51) \times10^{-4}$ & $(2.54\pm0.20) \times10^{-4}$ &
$(5.27\pm0.38) \times10^{-5}$ \\
& 0.4 & $(9.78\pm0.80) \times10^{-4}$ & $(2.60\pm0.20) \times10^{-4}$ &
$(6.58\pm0.67) \times10^{-5}$ \\
Direct & & $(8\pm3) \times10^{-4}$ & $(3\pm2) \times10^{-4}$ & 0
\\
\hline
\end{tabular*}
\end{table}
in the form of mean${}\pm{}$standard error, for the SISR [with several
choices of $\theta$ in
(\ref{uweights}), including $\theta=\wtd\theta_{2.5}$] and direct
Monte Carlo estimates of
$P_0 \{ S_n/n \geq2.5 \}$. It shows a variance reduction of 35 times for
$n=15$ and 80 times for $n=20$ over direct Monte Carlo when $\wtd
\theta_{2.5}$ is used as
an approximation to $\theta_{2.5}$ in the resampling weights (\ref
{uweights}) for SISR. When
$n=25$, direct Monte Carlo fails while the SISR estimate still has a
reasonably small
standard error.
\end{exa}

\begin{appendix}\label{app}
\section*{\texorpdfstring{Appendix: Proof of (\lowercase{\protect\ref{3.16}}) and (\lowercase{\protect\ref{ku}})}
{Appendix: Proof of (3.19) and (3.54)}}

\begin{pf*}{Proof of (\ref{3.16})}
For $0 < \varepsilon
< I$, let
\[
M_\varepsilon= \{ \theta\dvtx\phi(\mu_\theta)=I-\varepsilon\},\qquad
H(\theta) = \{ \mu\in\Lambda^o\dvtx\theta'(\mu-\mu_\theta) \geq0
\}.
\]
If $\mu\in H(\theta)$, then $\theta' \mu\geq\theta' \mu_\theta$
and, therefore,
%
%
\setcounter{equation}{0}
\begin{equation} \label{phibound}
\phi(\mu) = \sup_{\tilde\theta} \{ \wtd\theta' \mu- \psi(\wtd
\theta) \} \geq
\theta' \mu-\psi(\theta) \geq\theta' \mu_\theta- \psi(\theta)
= I-\varepsilon.
\end{equation}
Moreover, for $\theta\in M_\varepsilon$, $H(\theta)$ is a closed
half-space whose boundary is the
tangent space of $\{ \mu\dvtx\phi(\mu) = I-\varepsilon\}$
at $\mu_\theta$. Hence,
%
%
\begin{equation} \label{A.2}
\phi(\mu) \neq I-\varepsilon\qquad\mbox{for } \mu\in\Lambda^o \bigm\backslash
\bigcup_{\theta\in M_\varepsilon} H(\theta).
\end{equation}
Making use of this and (\ref{phibound}), we next show that
%
%
\begin{equation} \label{cupH}
\bigcup_{\theta\in M_\varepsilon} H(\theta) = \{ \mu\dvtx\phi(\mu)
\geq I-\varepsilon\}
\end{equation}
and, therefore, by (\ref{I}),
%
%
\begin{equation} \label{3.15}
\Gamma:= \{ \mu\dvtx g(\mu) \geq b \} \subset\{ \mu\dvtx\phi(\mu) \geq
I -\varepsilon\}
= \bigcup_{\theta\in M_\varepsilon} H(\theta).
\end{equation}
By (\ref{phibound}), $\bigcup_{\theta\in M_\varepsilon} H(\theta)
\subset\{ \mu\dvtx
\phi(\mu) \geq I-\varepsilon\}$. Therefore, it suffices for the proof
of (\ref{cupH}) to
show that $\{ \mu\dvtx\phi(\mu) < I-\varepsilon\} \supset\Lambda^o
\setminus
\bigcup_{\theta\in M_\varepsilon} H(\theta)$. Suppose this is not the case.
Then there exists $\mu_1 \in\Lambda^o \setminus\bigcup_{\theta\in
M_\varepsilon}
H(\theta)$ such that $\phi(\mu_1) \geq I-\varepsilon$. Since $\Lambda
^o \setminus
\bigcup_{\theta\in M_\varepsilon} H(\theta) \supset\{ \mu\dvtx\phi(\mu
) < I-\varepsilon\}$,
there exists $\mu_2 \in\Lambda^o \setminus\bigcup_{\theta\in
M_\varepsilon} H(\theta)$
such that $\phi(\mu_2) < I-\varepsilon$. By continuity of $\phi$,
there exists $\rho
\in(0,1)$ such that $\phi(\rho\mu_1+(1-\rho) \mu_2)=I-\varepsilon$. Since
$\Lambda^o \setminus H(\theta)$ is a half-space,
$\Lambda^o \setminus\bigcup_{\theta\in M_\varepsilon}
H(\theta) = \bigcap_{\theta\in M_\varepsilon}
(\Lambda^o \setminus H(\theta))$ is convex and, therefore, $\rho\mu
_1+(1-\rho) \mu_2 \in
\Lambda^o \setminus\bigcup_{\theta\in M_\varepsilon} H(\theta)$, but this
contradicts (\ref{A.2}), thereby proving (\ref{cupH}).

Define the measure $Q$ by
\[
\frac{dQ}{dP}({\mathbf Y}_n) = \int_M e^{\theta' S_n-n \psi(\theta)} \,d \theta/
\operatorname{vol}(M),
\]
where vol$(M)$ is the volume of $M$.
Let $\mu_n = S_n/n$ and $h_n(\theta) = \theta' \mu_n - \psi(\theta)$.
From (\ref{3.15}), it follows that if $\mu_n \in\Gamma$, then
there exists $\theta_* \in M_\varepsilon$ such that $\theta_*'(\mu_n -
\mu_{\theta_*})
\geq0$ and, therefore,
%
%
\begin{equation} \label{hstar}
h_n(\theta_*) = \theta_*' \mu_n - \psi(\theta_*) \geq\theta_*'
\mu_{\theta_*} -
\psi(\theta_*) = \phi(\mu_{\theta_*}) = I-\varepsilon,
\end{equation}
since $\theta_* \in M_\varepsilon$.
Let $B_n = \{ \theta\dvtx(\theta-\theta_*)' \nabla h_n(\theta_*) \geq
0, \| \theta- \theta_* \| \leq n^{-1/2} \}$.
Then for all $\theta\in B_n$,
$h_n(\theta) = h_n(\theta^*)+(\theta-\theta^*)' \nabla h_n(\theta_*)-
(\theta-\theta_*)' \nabla^2 \psi(\theta_*)(\theta-\theta_*)/2
+o(\| \theta-\theta_* \|^2)$ and, therefore, by (\ref{hstar}) and
the definition of $B_n$,
\[
h_n(\theta) \geq I - \varepsilon- (K+1)/(2n) \qquad\mbox{for all large } n,
\]
where $K = \sup_{\theta\in M} \| \nabla^2 \psi(\theta) \|$.
Hence, for all large $n$,
%
%
\begin{eqnarray} \label{volM}
\frac{dQ}{dP}(\bY_n) & \geq& {\mathbf1}_{\{ \mu_n \in\Gamma\}}
\int_{B_n} \exp\{ n h_n(\theta) \} \,d \theta/ \operatorname{vol}(M)
\nonumber\\[-8pt]\\[-8pt]
& \geq& {\mathbf1}_{\{ \mu_n \in\Gamma\}} (c_d/2) e^{nI-n \varepsilon
-(K+1)/2} n^{-d/2}/\operatorname{vol}(M), \nonumber
\end{eqnarray}
in which $c_d$ denotes the volume of the $d$-dimensional unit ball.
Letting $\varepsilon\rightarrow0$ in (\ref{volM}) yields
$(dQ/dP)(\bY_n) \geq e^{nI+o(n)} {\mathbf1}_{\{ \mu_n \in\Gamma\}}$
in which $o(n)$ is uniform in~$\bY_n$. Hence,
\begin{eqnarray*}
P \{ g(S_n/n) \geq b | \bY_k \} & = & E_Q \biggl[ \frac{dP}{dQ}
(\bY_n) {\mathbf1}_{\{ S_n/n \in\Gamma\}} \,\frac{dQ}{dP}(\bY_k) \Big| \bY_k
\biggr] \\
& \leq& e^{-n
I+o(n)} \,\frac{dQ}{dP}(\bY_k),
\end{eqnarray*}
proving (\ref{3.16}).
\end{pf*}

To prove (\ref{ku}), we use ideas similar to those in the
proof of Lemma 1 of~\cite{CL03}
and the following result of \cite{NN87}, page 568.
\begin{lem}\label{lemma3} Let $\tau(0)=0$. Under (\ref{m1}),
there exist regeneration times~$\tau(i)$, $i \geq1$, such that:

\begin{longlist}
\item
$\tau(i+1)-\tau(i)$, $i \geq0$, are i.i.d. random
variables,

\item$\{ X_{\tau(i)}, \ldots, X_{\tau(i+1)-1}, \xi
_{\tau(i)+1}, \ldots, \xi_{\tau(i+1)}
\}$, $i=0,1,\ldots,$ are independent blocks,

\item$X_{\tau(i)}$ has distribution $\nu$ for all
$i \geq1$.
\end{longlist}
\end{lem}
\begin{pf*}{Proof of (\ref{ku})}
Let $\wtd\ell_x
= E_x \{ \sum_{n=1}^{\tau} e^{\theta_b S_n} u(X_n) \}$, $\wtd\ell
_\nu= \int\wtd\ell_x
\,d \nu(x)$ and $A = \{ \tau(i)\dvtx i \geq1 \}$. Since $u \geq1$,
%
%
\begin{eqnarray} \label{expand}\quad
E_x \{ e^{\theta_b S_k} u(X_k) \} &=& E_x \bigl\{ e^{\theta_b
S_k} u(X_k) {\mathbf1}_{\{ \tau\geq k \}} \bigr\} \nonumber\\
&&{} + \sum_{j=1}^{k-1}
E_x\bigl(e^{\theta_b S_j} {\mathbf1}_{\{ j \in A \} }\bigr)
E_\nu\bigl(e^{\theta_b S_{k-j}} u(X_{k-j}) {\mathbf1}_{\{ \tau\geq
k-j \}}\bigr) \\
&\leq& \wtd\ell_x + \wtd\ell_\nu\Bigl[\sup_{j \geq1}
E_x\bigl(e^{\theta_b S_j}
{\mathbf1}_{\{ j \in A \}}\bigr)\Bigr]. \nonumber
\end{eqnarray}
Let $0 < \sigma=\sigma(1) < \sigma(2) < \cdots$ be the hitting
times of $C$. Then
%
%
\begin{eqnarray} \label{ltilde}
\wtd\ell_x &\leq& E_x \Biggl\{ \sum_{n=1}^{\sigma} e^{\theta_b S_n} u(X_n)
\Biggr\}\nonumber\\[-8pt]\\[-8pt]
&&{} + E_x \Biggl\{ \sum_{k\dvtx\sigma(k) < \tau} e^{\theta_b S_{\sigma(k)}}
\sum_{n=\sigma(k)+1}^{\sigma(k+1)} e^{\theta_b(S_n-S_{\sigma(k)})}
u(X_n) \Biggr\}.\nonumber
\end{eqnarray}
Let $y \in\cX$. By (U2), for all $n \geq2$,
\[
E_y \bigl\{ e^{\theta_b S_n} u(X_n) {\mathbf1}_{\{ n \leq\sigma\}} \bigr\}
\leq(1-\beta) E_y \bigl(e^{\theta_b S_{n-1}} u(X_{n-1}) {\mathbf1}_{\{ n-1
\leq
\sigma\}}\bigr),
\]
from which it follows by proceeding inductively and applying (U3) that
%
%
\begin{equation} \label{recursive}
E_y \Biggl\{ \sum_{n=1}^{\sigma} e^{\theta_b S_n} u(X_n) \Biggr\}
\leq\beta^{-1} \max\{ a, (1-\beta) u(y) \} \leq\alpha u(y),
\end{equation}
where $\alpha=\beta^{-1} \max\{ a, (1-\beta) \}$.
Substitution of (\ref{recursive}) into (\ref{ltilde}) then yields
%
%
\begin{equation} \label{lbound}
\wtd\ell_x\,{\leq}\,\alpha\Biggl\{\!u(x)\,{+}\,E_x \Biggl( \sum_{n=0}^{\tau-1}
e^{\theta_b S_n} u(X_n) {\mathbf1}_{\{ X_n \in C \}} \Biggr)\!\Biggr\}\,{\leq}\,\alpha u(x)
\,{+}\,\eta\ell_x(\theta_b;C),\hspace*{-45pt}
\end{equation}
where $\eta= \sup_{y \in C} u(y)$. Since $\int_\mathcal{X} u(x) \,d \nu(x) < \infty$
and $\int_\mathcal{X} \ell_x(\theta_b;C) \,d \nu(x) < \infty$, it
follows from (\ref{lbound})
that $\wtd\ell_\nu< \infty$. Combining
\[
\ell_x(\theta_b;C) \leq E_x (e^{\theta_b S_\sigma})
\Bigl[\sup_{y \in C} \ell_y(\theta_b;C)\Bigr]
\]
with (\ref{recursive}) yields
%
%
\begin{equation} \label{U3}
\sup_{x \in\cX} \{ \ell_x(\theta_b;C)/u(x) \} < \infty.
\end{equation}
Let $Q^*$ be a probability measure under which
\[
\frac{dQ^*}{dP_\nu}\bigl( \{ (X_t,S_t)\dvtx t \leq\tau(i) \}\bigr) = e^{\theta_b
S_{\tau(i)}}.
\]
Then
%
%
\begin{equation} \label{A.6}
\sup_{k \geq1} E_\nu\bigl(e^{\theta_b S_k} {\mathbf1}_{\{ k \in A \}}\bigr) =
\sup_{k \geq1}
Q^* \{ \tau(i)=k \mbox{ for some } i \} \leq1.
\end{equation}
From (\ref{expand}), (\ref{lbound}), (\ref{U3}) and
\begin{eqnarray*}
E_x \bigl(e^{\theta_b S_j} {\mathbf1}_{\{ j
\in A \}}\bigr) & = & E_x \bigl(e^{\theta_b S_\tau} {\mathbf1}_{\{ \tau=j \}}\bigr)
\\
& &{} + \sum_{h=1}^{j-1} E_x\bigl(e^{\theta_b S_\tau} {\mathbf1}_{\{
\tau=h \}}\bigr)
E_\nu\bigl(e^{\theta_b S_{j-h}} {\mathbf1}_{\{ j-h \in A \}}\bigr) \\
& \leq& E_x(e^{\theta_b S_\tau}) \Bigl\{ 1+\sup_{k \geq1} E_\nu
\bigl(e^{\theta_b
S_k} {\mathbf1}_{\{ k \in A \} }\bigr) \Bigr\},
\end{eqnarray*}
(\ref{ku}) follows from (\ref{A.6}).
\end{pf*}
\end{appendix}


%

%
\printaddresses

\end{document}